%% file: gkgrp.tex
\documentclass[10pt,oneside]{jm}
\usepackage{amsfonts,amssymb,amsmath,amscd,graphicx,rotating} 

\input{abbr}

\begin{document}
\raggedbottom
\pagestyle{myheadings}
\markboth{}{}
\thispagestyle{empty}
\setcounter{page}{1}
\begin{center} \Large\bf
On $p$-groups of Gorenstein-Kulkarni type
\vspace*{1em} \\
\large\rm
J\"urgen M\"uller and Siddhartha Sarkar
\vspace*{1em} \\
\it 
Dedicated to the memory of D. N. Verma
\vspace*{1em} \\
\end{center}

\begin{abstract} \noindent
A finite $p$-group is said to be of Gorenstein-Kulkarni type 
if the set of all elements of non-maximal order is a maximal subgroup. 
$2$-groups of Gorenstein-Kulkarni type arise naturally in the 
study of group actions on compact Riemann surfaces. In this paper, 
we proceed towards a classification of $p$-groups of Gorenstein-Kulkarni 
type. 
\end{abstract}

\tableofcontents

\input{intro}

\input{kulinv}

\input{trees}

\input{cycext}

\input{stems}

\input{2grpex}

\input{genex}


\abs\hrulefill

\abs\abs
{\sc J.M.:
Lehrstuhl D f\"ur Mathematik, RWTH Aachen \\
Templergraben 64, D-52062 Aachen, Germany} \\
{\sf Juergen.Mueller@math.rwth-aachen.de}

\abs
{\sc S.S.: 
Department of Mathematics, IISER Bhopal \\
ITI (Gas Rahat) Building, Govindpura \\ 
Bhopal 462\,023, Madhya Pradesh, India} \\
{\sf sidhu@iiserb.ac.in}

\end{document}

%% file: abbr.tex
\newcommand{\abs}{\vskip 0.5em\noindent}

\newcommand{\AbsT}[1]{\paragraph{\hspace{-1em} #1}}
\newcommand{\AbsTT}[1]{\paragraph*{#1}}

\newcommand{\Pf}{\paragraph*{Proof.}}

\newcommand{\ExT}[1]{\AbsT{Example: #1.}}

\newcommand{\Thm}{\AbsT{Theorem.}}
\newcommand{\ThmT}[1]{\AbsT{Theorem: #1.}}

\newcommand{\Cor}{\AbsT{Corollary.}}

\newcommand{\bfi}{\noindent {\bf i)} }

\newcommand{\bfii}{\noindent {\bf ii)} }

\newcommand{\bfiii}{\noindent {\bf iii)} }

\newcommand{\bfa}{\noindent {\bf a)} }
\newcommand{\bfb}{\noindent {\bf b)} }

\newcommand{\QED}{\hfill $\sharp$}

\newcommand{\N}{\mathbb N}
\newcommand{\Z}{\mathbb Z}
\newcommand{\Q}{\mathbb Q}

\newcommand{\cB}{\mathcal B}

\newcommand{\cK}{\mathcal K}

\newcommand{\tcO}{\widetilde{\mathcal O}}
\newcommand{\cO}{\mathcal O}

\newcommand{\cS}{\mathcal S}

\newcommand{\cT}{\mathcal T}

\newcommand{\cZ}{\mathcal Z}

\newcommand{\al}{\alpha}

\newcommand{\gm}{\gamma}

\newcommand{\dt}{\delta}

\newcommand{\ka}{\kappa}

\newcommand{\ph}{\varphi}

\newcommand{\Aut}{\text{Aut}}

\newcommand{\id}{\text{id}}

\newcommand{\Inn}{\text{Inn}}

\newcommand{\lcm}{\text{lcm}}

\newcommand{\Out}{\text{Out}}

\newcommand{\baraut}{\ov{\text{\hspace{0.5em}\rule{0em}{0.55em}}}}

\newcommand{\ld}{,\ldots\hskip0em ,}
\newcommand{\lr}[1]{\langle #1\rangle}
\newcommand{\ov}[1]{\overline{#1}}

\newcommand{\cl}[1]{\lceil #1\rceil}

\newcommand{\wti}[1]{\widetilde{#1}}
\newcommand{\wh}[1]{\widehat{#1}}
\newcommand{\mt}{\mapsto}

\newcommand{\ra}{\rightarrow}

\newcommand{\cn}{\colon}
\newcommand{\dcup}{\stackrel{.}{\cup}}

\newcommand{\spmid}{\,\mid\,}
\newcommand{\spnmid}{\,\nmid\,}

\newcommand{\bu}{\bullet}

\newcommand{\sseq}{\subseteq}
\newcommand{\smin}{\setminus}
\newcommand{\emp}{\emptyset}

\newcommand{\tm}{\times}

\newcommand{\GAP}{{\sf GAP}}

%% file: intro.tex
\section{Introduction}\label{intro}

\abs
The present work arose out of the study of actions
of finite $p$-groups on compact Riemann surfaces, and
how this is reflected in the structure of the groups acting.
Although the main focus of the present paper is on the
group theoretical side, we briefly recall the setting
in order to give some motivation for the developments presented later:

\AbsT{Genus spectra.}
Let $X$ be a compact Riemann surface of genus $g$,
and let $\Aut(X)$ be its automorphism group.
If $g\geq 2$, then by a famous theorem due to Hurwitz \cite{hur}
the group $\Aut(X)$ is finite, and has order $|\Aut(X)|\leq 84\cdot(g-1)$. 
More generally, a finite group $G$ is said to act on $X$ if $G$ can be 
embedded into $\Aut(X)$. Hence for any fixed Riemann surface $X$ 
of genus $g\geq 2$ Hurwitz's Theorem implies that there are only 
finitely many groups $G$ (up to isomorphism) acting on $X$.

\abs
But conversely, given a finite group $G$, there always is an infinite set 
$\text{spec}(G)$ of integers $g\geq 2$, called the {\bf genus spectrum} 
of $G$, such that there is a Riemann surface of genus $g$ 
being acted on by $G$, see \cite{kul,mm}. 
The problem of determining $\text{spec}(G)$ is called
the {\bf Hurwitz problem} associated with $G$, see \cite{mm},
where the minimum of $\text{spec}(G)$, being called
the {\bf strong symmetric genus} of $G$, is of particular interest.
For more details and the state of the art for various classes
of groups we refer the reader to \cite{Breuer,sar} 
and the further references in there.

\abs
To attack the Hurwitz problem, in \cite{kul} a group theoretic 
invariant $N(G)\in\N$, now called the {\bf Kulkarni invariant} of $G$, 
is introduced, such that 
\begin{equation*}\tag{$\ast$}
g\equiv 1\pmod{N(G)}\quad\text{whenever}\quad g\in\text{spec}(G).
\end{equation*}
Hence we may define the {\bf reduced genus spectrum} of $G$ as
$$ \text{spec}_0(G):=
   \left\{\frac{g-1}{N(G)}\in\N;g\in\text{spec}(G)\right\}\sseq\N .$$
Then it is also shown in \cite{kul} that the the complement 
$\N\smin\text{spec}_0(G)$ is finite, hence $N(G)$ appropriately 
describes the asymptotic behaviour of $\text{spec}(G)$. Moreover, this
also says that $N(G)$ is the (unique) maximal integer 
such that $(\ast)$ holds.


\AbsT{Groups of GK type.}
As was already said, the Kulkarni invariant $N(G)$ is
of group theoretic nature, where the details are given in
our restatement of Kulkarni's Theorem \cite{kul} in \ref{kulkarnithm}. 
The essential ingredient is a structural property 
of the Sylow $2$-subgroups of $G$.
As it turns out, the necessary notions can be defined for any
finite $p$-group, where $p$ is a rational prime:

\abs
Let $G$ be a finite $p$-group, and let 
$$ \exp(G):=\lcm\{|x|\in\N;x\in G\}=\max\{|x|\in\N;x\in G\} $$
be the exponent of $G$. Then $G$ is said to be of 
{\bf Gorenstein-Kulkarni type}, or of {\bf GK type} for short,
if the set
$$ \cK(G):=\{x\in G;|x|<\exp(G)\}\sseq G $$ 
of all elements of non-maximal order is a maximal subgroup of $G$. 
In this case, $\cK(G)\lhd G$ of course is a normal subgroup, 
and is called the {\bf GK kernel} of $G$.

\abs
In particular, any non-trivial cyclic $p$-group is of GK type,
but the trivial group is not.
We note that for $p=2$ this is
essentially the notion of groups of `type II' defined in \cite{kul},
except that there the cyclic groups were excluded.
Of course, the notion of groups of `type II' is the motivation for the 
considerations made here, and the name `Gorenstein-Kulkarni type' is
reminiscent of \cite{kul} and the acknowledgements made in there.

\AbsT{Other classes of $p$-groups.}
Groups of GK type will be the main objects of study in the
present paper.
But right now it seems to be worth-while to discuss briefly 
the relationship of the GK property to a few other, 
well-known group theoretical notions for $p$-groups: 

\abs\bfi
Let $G$ be a finite $p$-group.
Then $G$ is said to have the {\bf maximum exponent property (MEP)},
if the set $\cK(G)\sseq G$ is a subgroup (which then of course is normal)
such that $G/\cK(G)$ is abelian.
This notion was introduced in \cite{mt}, excluding the case $p=2$;
and of course any group of GK type has MEP.

\abs\bfii
For $i\in\N_0$ let 
$$ G^{p^i}:=\lr{x^{p^i}\in G;x\in G}\unlhd G .$$
Then $G$ is called {\bf powerful (POW)}, if either $p$ is odd 
and $G/G^p$ is abelian, or $p=2$ and $G/G^4$ is abelian,
see \cite[Ch.6.1]{lmc}. In particular, any abelian group is powerful,
hence this notion can be seen as a (proper) generalisation of 
being abelian. Moreover, any powerful group has MEP,
as is seen by essentially repeating the argument given 
in \cite{mt} for the case $p$ odd:

\abs
If $G$ is powerful, then by \cite[La.6.1.9]{lmc} 
$$ G^{p^i}/G^{p^{i+1}}\ra G^{p^{i+1}}/G^{p^{i+2}}\cn 
   xG^{p^{i+1}}\mt x^pG^{p^{i+2}} $$
is an epimorphism for all $i\in\N_0$.
Hence letting $\exp(G)=p^e$ for some $e\in\N_0$,
we conclude that $G\ra G\cn x\mt x^{p^{e-1}}$ is a homomorphism, 
whose kernel is $\cK(G)$, which hence is a subgroup.
Letting $G^{(1)}\unlhd G$ denote the derived subgroup of $G$,
we from $G^{(1)}\leq G^p\leq\cK(G)$ for $p$ odd,
and $G^{(1)}\leq G^4\leq G^2\leq\cK(G)$ for $p=2$,
infer that $G/\cK(G)$ is abelian, thus $G$ has MEP.
\QED

\abs\bfiii
The group $G$ is called {\bf regular (REG)}, if 
for all $x,y\in G$ we have
$$ (xy)^p=x^py^pz,\quad\text{for some }z\in((\lr{x,y})^{(1)})^p ,$$
see \cite[Sect.III.10]{hup}. In particular, any abelian group
is regular, hence this notion can also be seen as generalisation 
of being abelian; note that this is a proper generalisation 
for $p$ odd, while for $p=2$ by \cite[Thm.III.10.2]{hup} 
the regular groups are precisely the abelian ones.
Moreover, for $p$ odd neither of the notions of being regular
and being powerful implies the other; in particular
Wielandt's example reproduced in \cite[Sect.III.10.3]{hup} is
powerful but not regular.

\abs\abs
Thus, in view of the the implications just mentioned,
the above classes of $p$-groups are related to each other
as depicted in the following Venn diagram, where {\bf AB} 
denotes the class of abelian groups, and for $p=2$ 
the region {\bf REG} has to be deleted; but we point out
that we do not try to represent the various classes in any 
sense according to their (asymptotic) size:

\abs
\begin{center}
\includegraphics[height=40mm]{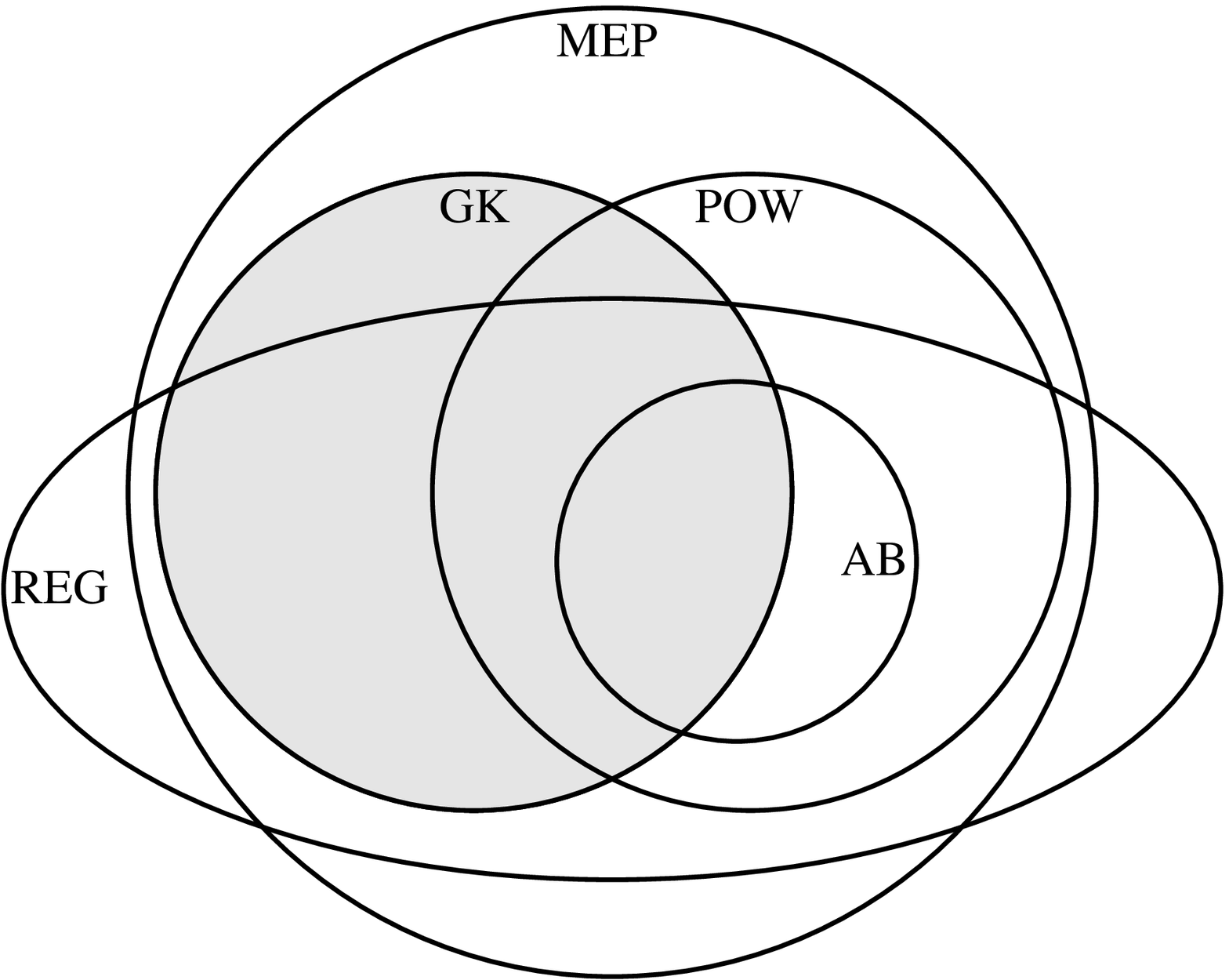}
\end{center}

\abs
Indeed, there are no further general implications between these properties,
that is all the regions depicted are actually non-empty,
as can be verified by the search techniques described in
Section \ref{2grpex}. We just mention the following examples:
The elementary abelian group $C_p^2$, 
since $\cK(C_p^2)=\{1\}$, is not of GK type;
by \ref{cp3tree} there are non-abelian groups of shape $C_p^3.C_p$
which are of GK type but are not powerful;
and the extraspecial group $E_+(p^{2+1})$ of order $p^3$ and
exponent $p$, where $p$ is odd, is regular, but since 
$\cK(E_+(p^{2+1}))=\{1\}$ does not have MEP.

\AbsT{GK trees.}
Hence groups of GK type, to all of our knowledge, 
form a new class of $p$-groups,
where it remains to be seen whether this indeed is an interesting one.
In particular, it is unclear whether this notion is of 
relevance outside the realm of $p$-groups; for example we are wondering
whether and how the structure of a finite group is influenced by
the fact that a Sylow ($2$-)subgroup is of GK type or not.

\abs
The purpose of the present paper is to understand the properties
of groups of GK type, and to set up some machinery to proceed towards their
classification. The basic idea is to organise the groups of GK type,
for a fixed rational prime $p$,
into a directed graph whose directed edges are given by connecting
a group of GK type to its GK kernel. It is immediate that its 
connected components are trees, where since walking along directed paths
amounts to iterating the step of taking a GK kernel, any such
GK tree is rooted at a group \textit{not} of GK type.

\abs
The GK trees considered here are, of course, modelled
after the so-called co-class trees used in the classification
business of \textit{all} $p$-groups; as general references 
for details about co-class theory and the structure 
of co-class trees see for example \cite{lmc} and \cite{EiLG}, respectively.
The most noticeable difference between the former and the latter, 
apart from a slightly changed terminology, is that our GK trees 
grow from bottom up, while co-class trees grow from top down.
This behaviour of GK trees (being much more sensible
from a biological viewpoint anyway) is due to the fact 
that GK groups are related to each other by forming (cyclic) 
upward extensions, while groups with fixed co-class are related 
to each other by forming (central) downward extensions.

\abs
Still, quite a few similarities remain, motivating the 
following questions; formal definitions of the relevant 
notions will be given in Section \ref{trees}:
$\bu$
Are there finite GK trees? 
$\bu$ 
Does a finite GK tree always consist of a single vertex?
$\bu$
Does an infinite GK tree have finitely many `stems' 
(called `main-lines' in \cite{EiLG})?
$\bu$
Does an infinite GK tree even have a single stem?
$\bu$
How do the `bushes' (called `branches' in \cite{EiLG}) 
of a GK tree look like? 
$\bu$
Is a GK tree always `periodic' (in the sense of \cite{EiLG})?

\abs
In the present paper we will answer the question of finiteness
affirmatively, and achieve a complete description of the stems 
of a GK tree. In particular, it will be shown that there are
finite GK trees having more than one vertex, and that an infinite
GK tree might have more than one, but always has finitely many stems.
A discussion of the bushes and of periodicity of GK trees is 
left to the sequel \cite{ms} to the present paper.

\AbsTT{Outline.}
The paper now is organised as follows: 
$\bu$
In Section \ref{kulinv}, we begin our journey by
reformulating Kulkarni's Theorem in terms of groups of GK type,
and prove a reduction argument. Actually, this was the original 
incentive for the present work, and already contains the germ 
of relating groups of GK type to their iterated GK kernels.
We conclude by giving a simple statement
relating the group structure of a finite group to the
GK property of its Sylow $2$-subgroups.
$\bu$
In Section \ref{trees} 
we introduce groups of GK type formally, together 
with some associated group theoretic notions, and collect
a few immediate properties. Moreover, we introduce 
GK trees formally, together with the relevant
graph theoretic notions.
$\bu$
In Section \ref{cycext} we develop the necessary machinery
of cyclic extensions.
$\bu$
In Section \ref{stems} we give a group theoretical 
characterisation of whether a $p$-group not of GK type 
is the root of an infinite GK tree. Moreover, in this case,
we describe the branching behaviour of the stems,
and the groups attached to the vertices lying on the stems.
$\bu$
In Section \ref{2grpex} we present a collection of
examples consisting of $2$-groups. These largely have been found in the 
first place by searching the {\sf SmallGroups} database \cite{SmGrp} 
available in through the computer algebra system \GAP{} \cite{GAP}. 
$\bu$
In Section \ref{genex}, finally, we present a collection of
`generic' examples, in the sense that the prime $p$ is treated as
a parameter.

\abs
The reader is recommended to keep the trees depicted in 
Tables \ref{treetbl1}--\ref{treetbl10} in mind while 
going through the theoretical parts 
of the paper. For the relevant background from general group theory 
and the theory of $p$-groups, we refer for example to 
\cite{hup} and \cite{lmc}, respectively.

\AbsTT{Acknowledgements.}
We both gratefully acknowledge generous hospitality 
of the Harish-Chandra Research Institute (HRI), Allahabad,
where we had the opportunity to meet D. N. Verma, 
discussions with whom left us deeply impressed.
Moreover, important parts were carried out when 
the second author was enjoying hospitality of
Professor Aner Shalev at Hebrew University Jerusalem, supported
by a Golda Meir Fellowship (2008--2009), and hospitality of
Professor Gerhard Hiss at RWTH Aachen (September 2008).
Finally, we would like to thank Avionam Mann for helpful suggestions.

%% file: kulinv.tex
\section{The Kulkarni invariant}\label{kulinv}

\abs
In this section we reformulate Kulkarni's Theorem 
using the language of groups of GK type. 
We then proceed to give a reduction, based on group theoretic transfer,
of the proof of Kulkarni's Theorem to the case where the Sylow 
$2$-subgroups of the group under consideration are \textit{not} of GK type.
Actually, proving this reduction was the initial ignition to pursue
the idea of relating $p$-groups of GK type to each other by going over to GK 
kernels and iterating this step; we comment on this at the end of \ref{reduc}.

\abs
Unfortunately, in the case of Sylow $2$-subgroups not of GK type
pure group theory does not seem to provide a conceptually better proof 
of Kulkarni's Theorem than the one using the combinatorics of the 
Riemann-Hurwitz equation already presented in \cite{kul}. 
We just note that, since non-trivial cyclic groups are of GK type 
but are not of `type II', our approach still leads to an improvement 
by allowing to avoid the case distinction in \cite[Sect.2.10]{kul}.

\abs
Anyway, the Riemann-Hurwitz equation is the 
key tool to relate the property of $G$ acting on a Riemann surface 
of genus $g\geq 2$ to a group theoretic invariant; we recall
the relevant theorem, again due to Hurwitz \cite{hur}, for convenience:

\ThmT{(Hurwitz \cite{hur})}\label{rhthm}
Let $G$ be a finite group. Then
$G$ acts on a Riemann surface of genus $g\geq 2$ 
if and only if there are $h,r\in\N_0$ and integers $n_1\ld n_r\geq 2$, 
fulfilling the {\bf Riemann-Hurwitz equation}
$$ 2\cdot(g-1)
=|G|\cdot\left(2\cdot(h-1)+\sum_{i=1}^r(1-\frac{1}{n_i})\right) ,$$
such that there exits a generating set $\{a_1,b_1\ld a_h,b_h;c_1\ld c_r\}$
of $G$ such that $|c_i|=n_i$, for all $i\in\{1\ld r\}$, and 
fulfilling the {\bf long relation} 
$$ [a_1, b_1]\cdot\cdots\cdot [a_h, b_h]\cdot c_1\cdot\cdots\cdot c_r=1 ;$$
here $[x,y]:=x^{-1}y^{-1}xy\in G$ denotes the commutator of $x,y\in G$.

\abs
The elements of $\{a_1,b_1\ld a_h,b_h\}$ and $\{c_1\ld c_r\}$ are 
called {\bf hyperbolic} and {\bf elliptic} generators, respectively.

\ThmT{(Kulkarni \cite{kul})}\label{kulkarnithm} 
Let $G$ be a finite group, and for any rational prime $p$ 
let $G_p$ denote a Sylow $p$-subgroup of $G$. Moreover, let
$$ \gm(G):=\left\{\begin{array}{rl}
2,&\text{if }G_2\neq\{1\}\text{ and }G_2\text{ is not of GK type},\\
1,&\text{if }G_2=\{1\}\text{ or }G_2\text{ is of GK type}.\\
\end{array}\right. $$
Then the Kulkarni invariant of $G$ is given as 
$$ N(G)=\frac{1}{\gm(G)}\cdot\prod_{p\spmid |G|}\frac{|G_p|}{\exp(G_p)} .$$

\AbsT{Reduction to the case where $G_2$ is not of GK type.}\label{reduc}
We present a reduction of the proof of \ref{kulkarnithm} to the
case where $G_2$ is not of GK type. Hence we may assume the assertion 
to be true if $G_2$ is not of GK type. The approach is based
on group theoretic transfer, which we assume the reader to be 
familiar with; as a general reference see for example \cite[Sect.IV]{hup}.

\abs\bfi
Let $G_2$ be of GK type, let 
$\cK(G_2):=\{x\in G_2;|x|<\exp(G_2)\}\lhd G_2$ be its GK kernel,
let $\pi\cn G_2\ra G_2/\cK(G_2)\cong C_2\cong\{\pm 1\}$
be the natural epimorphism, and let $\ph\cn G\ra\{\pm 1\}$ 
be the associated transfer homomorphism with kernel $H:=\ker(\ph)\unlhd G$.
Recall that any element $x\in G$ can be uniquely written as $x=x_2 x_{2'}$,
where $x_2,x_{2'}\in G$ such that $|x_2|$ is a $2$-power, 
$|x_{2'}|$ is odd and $x_2 x_{2'}=x_{2'} x_2$.
An element $x$ is said to have {\bf $2$-maximal order} if 
$|x_2|=\exp(G_2)=\exp_p(G)$, the $2$-exponent of $G$.
We first show that 
$$ H=\{x\in G;x\text{ does not have $2$-maximal order}\} :$$

\abs
In order to do so, let $x\in G$.
Since $|x_{2'}|$ is odd, we have $\ph(x_{2'})=1$, implying that
$\ph(x)=\ph(x_2)$, thus to determine $\ph(x)\in\{\pm 1\}$
we may assume that $x=x_2$ is a $2$-element.
Moreover, there are $y_1\ld y_t\in G$ and $r_1\ld r_t\in\N$, 
for some $t\in\N$, such that $\sum_{i=1}^t r_i=[G\cn G_2]$ and 
$$ y_i x^{r_i} y_i^{-1}\in G_2,\quad\text{for all }i\in\{1\ld t\},
\quad\text{and}\quad\ph(x)=\prod_{i=1}^t\pi(y_i x^{r_i} y_i^{-1}) .$$
Since $|y_i x^{r_i} y_i^{-1}|=|x^{r_i}|$, for all $i\in\{1\ld t\}$, 
we conclude that $y_i x^{r_i} y_i^{-1}\in\cK(G_2)$
whenever $x$ does not have $2$-maximal order, and thus we have
$\ph(x)=1$ in this case. If $x$ has $2$-maximal order, then we similarly 
have $y_i x^{r_i} y_i^{-1}\in\cK(G_2)$ if and only if $r_i$ is even,
thus letting 
$$ l:=|\{i\in\{1\ld t\};r_i\text{ is odd}\}|\in\N_0 $$ 
we get $\ph(x)=(-1)^l$, where from $\sum_{i=1}^t r_i=[G\cn G_2]$ being odd 
we conclude that the cardinality $l$ is odd as well, 
implying that $\ph(x)=-1$ in this case.
This shows that $H$ is as asserted;
in particular we have $H\lhd G$ such that $G/H\cong C_2$.

\abs\bfii
Now, for $g\in\text{spec}(G)$ we have to show that 
$g\equiv 1\pmod{N(G)}$, where $N(G)\in\N$ is given by the 
formula in the assertion:
To this end, we choose a generating set $\{a_1,b_1\ld a_h,b_h;c_1\ld c_r\}$ 
of $G$ as in Hurwitz's Theorem \ref{rhthm}. Still writing
$n_i:=|c_i|$ for all $i\in\{1\ld r\}$, dividing both sides 
of the associated Riemann-Hurwitz equation by $2\cdot N(G)$ yields
$$ \frac{g-1}{N(G)}=\left(\prod_{p\spmid |G|}\exp(G_p)\right)\cdot
   \left((h-1)+\sum_{i=1}^r\frac{n_i-1}{2n_i}\right)\in\Q .$$
Thus we have to show that the latter expression actually is an integer:

\abs
The $i$-th summand 
$\frac{n_i-1}{2n_i}\cdot\prod_{p\spmid |G|}\exp(G_p)\in\Q$
is an integer if $c_i$ does not have $2$-maximal order,
while if $c_i$ has $2$-maximal order, then $n_i$ is even and thus
the number $\frac{n_i-1}{2n_i}\cdot\prod_{p\spmid |G|}\exp(G_p)\in\Q$
is half an integer but not an integer. Hence we have to show that
the cardinality
$$ m:=|\{i\in\{1\ld r\};c_i\text{ has $2$-maximal order}\}|
     =|\{i\in\{1\ld r\};\ph(c_i)=-1\}| $$ 
is even:
Applying the transfer homomorphism $\ph\cn G\ra\{\pm 1\}$ to the
long relation associated with the generating set 
$\{a_1,b_1\ld a_h,b_h;c_1\ld c_r\}$,
and using $G^{(1)}\leq H$, we get $\ph(c_1)\cdot\cdots\cdot\ph(c_r)=1$,
implying that $m$ is even.

\abs\bfiii
We finally proceed to show that for all but finitely many $g\geq 2$
satisfying $g\equiv 1\pmod{N(G)}$ we actually have $g\in\text{spec}(G)$:
In order to do so, we may assume that the assertion has been proved for $H$, 
since if a Sylow $2$-subgroup $H_2$ of $H$ is of GK type again
we may proceed by induction, while if $H_2$ is not of GK type we 
assume the assertion to be true anyway; note that in particular 
the trivial group is not of GK type, but the assertion holds 
for that group anyway.

\abs
Hence by assumption there is $k\in\N$ such that we have 
$g'\in\text{spec}(H)$ whenever $g'\geq 2$ satisfies
$$ g'\equiv 1\pmod{N(H)}\quad\text{and}\quad\frac{g'-1}{N(H)}\geq k ,$$
where $N(H)$ is given by the formula in the assertion.
For any such $g'$ there is a generating set 
$\{a_1,b_1\ld a_h,b_h;c_1\ld c_r\}$ of $H$
as in Hurwitz's Theorem \ref{rhthm}.
Extending the generating set of $H$ 
by hyperbolic generators $a_{h+1}=b_{h+1}\in G\smin H$,
we again obtain a generating set of $G$ fulfilling the long relation,
and the Riemann-Hurwitz equation yields
$$ 2\cdot(g-1)=|G|\cdot\left(2h+\sum_{i=1}^r(1-\frac{1}{n_i})\right)
\!=2\cdot|H|\cdot\left(2 +\frac{2\cdot(g'-1)}{|H|}\right)
=4\cdot(|H|+g'-1) .$$
From $\exp(G_p)=\exp(H_p)$ for $p$ odd, and $\exp(G_2)=2\exp(H_2)$, we get
$N(G)=\gm(G)\cdot N(G)=\gm(H)\cdot N(H)$.
This implies 
$$ \frac{g-1}{N(G)}=\frac{2\cdot(|H|+g'-1)}{\gm(H)\cdot N(H)} 
=\frac{2}{\gm(H)}\cdot\frac{g'-1}{N(H)}+\prod_{p\spmid |G|}\exp(G_p) ,$$ 
where the second summand is even.

\abs
Extending the generating set of $H$ by elliptic generators 
$c_{r+1}=c_{r+2}^{-1}\in G\smin H$ instead,
we obtain a generating set of $G$ fulfilling the long relation,
where the Riemann-Hurwitz equation this time yields
$$ 2\cdot(g-1)
=2\cdot|H|\cdot\left(2\cdot\frac{n_{r+1}-1}{n_{r+1}} 
+\frac{2\cdot(g'-1)}{|H|}\right)
=4\cdot\left(\frac{n_{r+1}-1}{n_{r+1}}\cdot |H|+g'-1\right) ,$$
implying
$$\frac{g-1}{N(G)}=\frac{2}{\gm(H)}\cdot\frac{g'-1}{N(H)}+
 \frac{n_{r+1}-1}{n_{r+1}}\cdot\prod_{p\spmid |G|}\exp(G_p) ,$$ 
where since $c_{r+1}$ has $2$-maximal order the second summand is odd.

\abs
In conclusion, since $\frac{2}{\gm(H)}\in\{1,2\}$, we have 
$g\in\text{spec}(G)$ whenever $g\geq 2$ satisfies
$$ \rule{3.9em}{0em}
g\equiv 1\pmod{N(G)}\quad\text{and}\quad
\frac{g-1}{N(G)}\geq\frac{2}{\gm(H)}\cdot k+\prod_{p\spmid |G|}\exp(G_p) .
\rule{3.5em}{0em} \sharp $$

\abs\abs
We remark that, in part (iii) of the above proof, 
from $G_2/(G_2\cap H)\cong G_2H/H = G/H \cong C_2$
we infer that the Sylow $2$-subgroup of $H$ can be
chosen as 
$$ H_2:=G_2\cap H 
   =\{x\in G_2;x\text{ does not have $2$-maximal order}\}=\cK(G_2) .$$
Hence a single step in the above reduction process
amounts to taking a GK kernel, 
and the induction argument just says to iterate this.

\Cor\label{kulkarnicor}
Let $G$ be a $2$-perfect group, that is $[G\cn G^{(1)}]$ is odd.
Then the Sylow $2$-subgroups of $G$ are not of GK type.

\Pf
Assume to the contrary that the Sylow $2$-subgroups of $G$
are of GK type. Then 
$H=\{x\in G;x\text{ does not have $2$-maximal order}\}$
is a normal subgroup of index $2$, a contradiction.
\QED

\abs\abs
We remark that a special case of \ref{kulkarnicor} is already
contained in \cite{kul}: Using the combinatorics of the 
Riemann-Hurwitz equation it is shown there that the Sylow $2$-subgroups 
of a perfect group $G$ are not of GK type.
Finally we note that the statement of \ref{kulkarnicor}
does not hold for $p$ odd: For example, the symmetric group $\cS_3$ 
is $3$-perfect, but its Sylow $3$-subgroups are cyclic, thus are of GK type.

%% file: trees.tex
\section{Groups of Gorenstein-Kulkarni type}\label{trees}

\abs
Having the new notion of groups of GK type in our hands,
we now set out to develop some theory to understand their
structure and to proceed towards a classification.
We begin by recalling the basic definition: 

\AbsT{Groups of GK type.}
\bfa
Let $p$ be a rational prime. 
A finite $p$-group $G$ is said to be of {\bf Gorenstein-Kulkarni type},
or {\bf GK type} for short, if the set 
$$ \cK(G):=\{g\in G;|g|<\exp(G)\}\sseq G $$ 
of all elements of non-maximal order is a maximal subgroup of $G$.

\abs
In this case, $\cK(G)\lhd G$ is a characteristic subgroup of $G$ 
of index $p$, being called the {\bf GK kernel} of $G$.
Hence $G$ is of GK type if and only if there is an 
epimorphism $\pi\cn G\ra C_p$ such that $\ker(\pi)=\cK(G)$.

\abs
Note that the condition on $\cK(G)$ might fail in various ways:
The set $\cK(G)$ might fail to be a subgroup, or $\cK(G)$ might
be a subgroup but fails to be maximal; for examples see \ref{2grpexmpl}.
In particular, the trivial group is not of GK type.

\abs\bfb
As was already mentioned earlier, the key tool to describe 
groups of GK type is the following idea:
If $G$ is of GK type, then its GK kernel $\cK(G)$ might be of 
GK type again. In this case we may iterate the process of 
taking GK kernels, until we end up with a group not being of GK type. 
More formally, letting 
$$\cK^0(G):=G\quad\text{and}\quad
\cK^{i+1}(G):=\cK(\cK^i(G)),\quad\text{for }i\in\N_0 ,$$
yields a strictly descending chain of subgroups, called the
{\bf GK series} of $G$,
$$ R:=\cK^l(G)<\cK^{l-1}(G)<\cdots<\cK^2(G)<\cK^1(G)=\cK(G)<\cK^0(G)=G ,$$
for some $l\in\N$, where $R$ is not of GK type.
Then $l=l(G)\in\N$ is called the {\bf GK level}, and 
$R$ is called the {\bf GK root} of $G$;
for completeness $R$ is given GK level $l(R):=0$. 
Moreover, $\cK^i(G)\unlhd G$ is a characteristic subgroup, for all 
$i\in\{0\ld l\}$, and $G$ is called a {\bf GK extension} of $\cK^i(G)$.

\abs
Note that, if $G$ has GK level $l=1$, then it might still 
have a maximal subgroup being of GK type, or 
if $G$ has GK level $l\geq 2$, then it might have maximal
subgroups different from $\cK(G)$ being of GK type; 
for examples see \ref{2grpexmplcont}.

\AbsT{Immediate properties of groups of GK type.}\label{immediate}
Let $G$ be a group of GK type of level $l\in\N$, having root $R:=\cK^l(G)$.

\abs\bfa
For any $x\in G\smin\cK(G)$ and $i\in\{0\ld l-1\}$ we have 
$x^{p^i}\in\cK^i(G)\smin\cK^{i+1}(G)$. This for $i\in\{0\ld l\}$
implies $\cK^i(G)=\lr{x^{p^i},R}$ and $G=\lr{x,\cK^i(G)}$; 
the latter can be rephrased as $G/\cK^i(G)\cong C_{p^i}$, 
that is $G$ is a cyclic extension of $\cK^i(G)$. 

\abs
Moreover, from
$[G\cn\cK^i(G)]=p^i$ and $\exp(\cK^i(G))=\frac{\exp(G)}{p^i}$ we infer 
$$\frac{|G|}{\exp(G)}=\frac{|\cK^i(G)|}{\exp(\cK^i(G))}
 =\frac{|R|}{\exp(R)}, 
\quad\text{for all }i\in\{0\ld l\} .$$
Hence we have $\dt(G)=\dt(\cK^i(G))=\dt(R)$, where $\dt(G)\in\N_0$ denotes
the {\bf cyclic deficiency} of $G$ being defined as
$$ \dt(G):=\log_p(\frac{|G|}{\exp(G)})=\log_p(|G|)-\log_p(\exp(G)) ;$$
note that we have $\dt(G)=0$ if and only if $G$ is cyclic.

\abs\bfb
For any finite $p$-group $H$ let $\Phi(H)\unlhd H$ be its 
Frattini subgroup. Recall that by Burnside's Basis Theorem
\cite[Thm.III.3.15]{hup} we have $H/\Phi(H)\cong C_p^{r(H)}$,
where $r(H)\in\N_0$ is called the {\bf (generator) rank} of $H$,
and coincides with the cardinality of any minimal generating set of $H$. 

\abs
Then we have $\Phi(G)\leq\cK(G)$, 
in particular implying that $\exp(\Phi(G))<\exp(G)$.
But we have $\Phi(G)\not\leq\cK^2(G)$ if $l\geq 2$,
implying that in this case $\Phi(\cK(G))\leq\cK^2(G)\cap\Phi(G)<\Phi(G)$. 
Moreover, we have
$$ 1\leq r(G)\leq r(\cK(G))\leq\cdots\leq r(\cK^{l-1}(G))\leq r(R)+1 :$$

\abs
Since $\cK^{l-1}(G)$ is a cyclic 
extension of $R$, we have $r(\cK^{l-1}(G))\leq r(R)+1$.
Now consider $\cK^i(G)$ for $i\in\{0\ld l-2\}$. 
Then $\cK^{i+1}(G)$ is of GK type again,
and for $g\in\cK^i(G)\smin\cK^{i+1}(G)$ we have 
$g^p\in\cK^{i+1}(G)\smin\cK^{i+2}(G)$.
Since $\cK^{i+2}(G)<\cK^{i+1}(G)$ is a maximal subgroup, we have
$g^p\not\in\Phi(\cK^{i+1}(G))$. Hence there is a minimal generating set of 
$\cK^{i+1}(G)$ containing $g^p$, 
and thus $r(\cK^i(G))\leq r(\cK^{i+1}(G))$.
\QED

\ExT{Abelian groups}\label{abelian}
The first examples we are tempted to look at are the abelian groups:
Let $G$ be an abelian $p$-group with abelian invariants
$(p^{e_1}\ld p^{e_r})$, where $r\in\N_0$ and
$0<e_1\leq\cdots\leq e_r$; we allow for $r=0$,
letting $e_0:=0$, to catch the case $G=\{1\}$.

\abs
If $G\neq\{1\}$, that is $r\geq 1$, then from $\exp(G)=p^{e_r}$
we conclude that the set $\cK(G)=\{x\in G;|x|<p^{e_r}\}$
is a subgroup of $G$, having abelian invariants
$(p^{e_1}\ld p^{e_{s-1}},p^{e_r-1}\ld p^{e_r-1})$,
where $s\in\{1\ld r\}$ is the unique number such that
$e_{s-1}<e_s=\cdots=e_r$.
Thus $G$ is of GK type if and only if $s=r\geq 1$,
or equivalently $r\geq 1$ and $e_{r-1}<e_r$.
In particular, for $r=1$ we recover the (obvious) fact that
any non-trivial cyclic $p$-group is of GK type.

\abs
If $G$ is of GK type, then $G$ has GK level $l:=e_r-e_{r-1}\in\N$,
and for all $i\in\{0\ld l\}$
the abelian invariants of $\cK^i(G)$ are
$(p^{e_1}\ld p^{e_{r-1}},p^{e_r-i})$.
In particular, the GK root $R:=\cK^l(G)$ of $G$ has abelian invariants
$(p^{e_1}\ld p^{e_{r-1}},p^{e_{r-1}})$ whenever $r\geq 2$,
while for $r=1$, that is $G\neq\{1\}$ is cyclic,
we get the trivial group $\{1\}$ as GK root.
In particular, we have $r(R)=r$ whenever $G$ has rank
$r(G)=r\geq 2$, but $r(R)=0$ if $r(G)=r=1$.

\AbsT{Trees and stems.}
Our overall aim now is to get an overview over the (isomorphism types of)
groups of GK type in terms of their roots: 

\abs\bfa 
To this end, 
for any (isomorphism type of) finite $p$-group $R$ not of GK type
we define a connected directed tree $\cT(R)$ rooted in $R$,
being called the associated {\bf GK tree}, as follows: 
The vertices of $\cT(R)$ are the root $R$, 
together with the (isomorphism types of) finite $p$-groups $G$ 
being of GK type and having $R$ as their root, 
and from any vertex $G$ of GK type of $\cT(R)$ 
precisely one edge emanates and this edge ends in $\cK(G)$.

\abs
An extended collection of explicit examples is given in Sections 
\ref{2grpex} and \ref{genex}, supported by the pictures in Tables 
\ref{treetbl1}--\ref{treetbl10}. These examples in particular show that
$\cT(R)$ might have precisely one vertex, namely only the 
root $R$, which just means that $R$ does not occur as a GK root,
and that there are GK trees having more than one, 
but finitely many vertices. Finally, except the obvious cases
in Table \ref{treetbl3}, the depicted finite segments seem to indicate
periodic behaviour, at least from some level on, hence in particular
these trees should be infinite; it will follow from \ref{trivextthm}
that the trees shown are indeed infinite.

\abs
We remark that, by \ref{immediate},
all the groups belonging to $\cT(R)$ have the same
cyclic deficiency, namely $\dt(R)$, which hence is an invariant associated 
to $\cT(R)$. But the rank of the groups occurring might indeed vary, 
following the general pattern described in \ref{immediate},
and more specially the pattern in \ref{trivext} in the 
`trivial' extension case discussed below;
we will have an eye on this in the explicit examples given later.

\abs\bfb
In order to describe the structure of infinite GK trees,
we first need the following notion: Note first that, since 
$\cT(R)$ is a tree, for any $G$ in $\cT(R)$ there is a unique 
directed path in $\cT(R)$, say
$$ G=G_l\ra G_{l-1}\ra\cdots\ra G_1\ra G_0=R ,$$
connecting $G$ with the root $R$;
its length $l\in\N_0$ coincides with the GK level of $G$.
Moreover, $\cT(R)$ is {\bf locally finite} in the sense that 
there are only finitely many directed paths of any fixed length,
or equivalently there only finitely many 
(isomorphism types of) groups in $\cT(R)$ of any fixed GK level.

\abs
An infinite directed path in $\cT(R)$, say
$$ \cdots\ra G_2\ra G_1\ra G_0=R ,$$ 
ending in the root $R$ is called a {\bf stem} of $\cT(R)$.
Then, since $\cT(R)$ is locally finite, we conclude that $\cT(R)$ is infinite 
if and only if it has directed paths of arbitrarily large length, which
in turn is equivalent to $\cT(R)$ having a stem.

\ExT{The trivial group}\label{trivgrp}
The prototypical example of an infinite GK tree is the tree $\cT(\{1\})$
rooted at the trivial group: As was already discussed in \ref{abelian},
if $G\neq\{1\}$ is cyclic then $G$ belongs to $\cT(\{1\})$,
and since conversely any GK extension of the trivial group is cyclic, 
the GK extensions of the trivial group are precisely
the cyclic groups $C_{p^l}$, for all $l\in\N$, where $l$ coincides
with the GK level.
Thus $\cT(\{1\})$ is infinite, but as such is as simple as possible, 
just being a single stem without any branching, just consisting
of the directed edges $C_{p^l}\ra C_{p^{l-1}}$;
see Table \ref{treetbl1}.

\AbsT{Bushes.}
Let $R$ be a finite $p$-group not of GK type, and let $\cT(R)$
be the associated GK tree. In order to describe the vertices of $\cT(R)$ 
not lying on any of the stems we introduce new invariants as follows: 

\abs\bfa
Let $G$ in $\cT(R)$ be connected to the root $R$ by the directed path
$$ G=G_l\ra G_{l-1}\ra\cdots\ra G_1\ra G_0=R $$ 
of length $l\in\N_0$.
Then we let the {\bf height} $s=s(G)\in\{0\ld l\}$ of $G$ 
be the largest number such that $G_s$ lies on a stem of $\cT(R)$.
Moreover, the {\bf depth} $t=t(G):=l-s\in\N_0$ of $G$
coincides with the distance of $G$ from any of the stems of $\cT(R)$;
in particular, $G$ has depth $t=0$ if and only if it lies on a stem.

\abs
Next, for any vertex $H$ in $\cT(R)$, 
let $\cT(H)$ be the full sub-tree of $\cT(R)$ rooted in $H$, 
that is the tree consisting of the vertices in $\cT(R)$
possessing a directed path to $H$, and the directed edges of $\cT(R)$ 
between them; hence $\cT(H)$ is infinite if and only if $H$ lies
on a stem of $\cT(R)$.

\abs
Now, let the {\bf bush} $\cB(H)$ rooted in $H$ be the sub-tree 
of $\cT(H)$, in turn, consisting of the vertex $H$ together with 
the vertices of $\cT(H)$ not lying on any of the stems of $\cT(R)$,
and the directed edges between them. 
Thus $\cB(H)$ is finite, where we more formally have
$$ \cB(H):=\{H\}\cup\{G\in\cT(H);t(G)\geq 1\} .$$
Hence for the most interesting case of $H$ 
lying on a stem of $\cT(R)$ we have
$$ \cB(H)=\{H\}\dcup\{G\in\cT(R);G_{s(G)}=H\} ,$$
that is, in prosaic words, $\cB(H)$ captures precisely
the vertices $G$ of $\cT(R)$ being reached by branching off
a stem of $\cT(R)$ at the vertex $H$, and the depth of $G$
coincides with the distance of $G$ from the root $H$ of $\cB(H)$.

\abs
For completeness, if $\cT(R)$ is finite, then for any vertex
$G$ in $\cT(R)$ we let $s=s(G):=0$ and $t=t(G):=l\in\N_0$;
in particular, $G$ has depth $t=0$ if and only if $G=R$.
Moreover, for the most interesting case of the bush rooted at $R$
we just recover $\cB(R)=\cT(R)$, that is $\cT(R)$ is just the 
bush rooted in $R$.

\abs\bfb
Finally, we turn to periodicity:
Let $\cT(R)$ be infinite, such that there is precisely one stem,
$$ \cdots\ra H_2\ra H_1\ra H_0=R .$$
Given $k\in\N$, the tree $\cT(R)$ is called {\bf $k$-periodic} 
from level $l\in\N_0$ on, or just {\bf periodic} in the case $k=1$,
if for all $s\geq l$ the bushes 
$\cB(H_s)$ and $\cB(H_{s+k})$ are isomorphic as directed trees.

\abs
We remark that, if $\cT(R)$ has more than one but finitely many stems,
we may go over to suitable full sub-trees $\cT(H)$ instead, where $H$ lies
on a stem of $\cT(R)$ and has a sufficiently large level, so that
$\cT(H)$ has only one stem; actually, as will be shown in \ref{stemrem},
any infinite GK tree indeed has only finitely many stems.

\abs\abs
We are now prepared to specify the further programme:
The overall aim, of course, is to analyse the structure of GK trees.
In the present paper we proceed to give a group theoretic criterion 
on the root group deciding whether the associated GK tree is infinite,
and in this case describing the branching behaviour of the stems, 
and the groups attached to the vertices belonging to the stems.
In the sequel \cite{ms} of the present paper we will 
then tackle the description of the bushes, 
and deal with questions of periodicity.

%% file: cycext.tex
\section{Cyclic extensions}\label{cycext}

\abs
We now proceed to develop a framework to describe GK extensions.
To do so, we first, without further assumptions, look at cyclic 
extensions in general.

\AbsT{Cyclic extensions.}\label{cycextdef}
Let $H$ be a finite group, and
let $G$ be a cyclic extension of $H$ of degree $d\in\N$,
that is we have $H\lhd G$ such that $G/H\cong C_d$;
the cyclic extension is called {\bf proper} if $d>1$.
Let $g\in G$ such that $G=\lr{g,H}$; this is equivalent
to saying $\lr{\ov{g}}=G/H$,
where $\baraut\cn G\ra G/H$ denotes the natural epimorphism.

\abs
Let $\ka_g\in\Aut(H)$ be the conjugation automorphism
$$ \ka_g\cn H\ra H\cn y\mt y^g:=g^{-1}yg $$
induced by $g$.
Hence letting $h:=g^d\in H$ we infer that 
$(\ka_g)^d=\ka_h\in\Inn(H)$ is an inner automorphism,
and thus $\ov{\ka}_g\in\Out(H):=\Aut(H)/\Inn(H)$ 
has order dividing $d$,
where again $\baraut\cn\Aut(H)\ra\Out(H)$ denotes the natural epimorphism.
 
\abs
The cyclic extension is called {\bf trivial} if 
$\ov{\ka}_g=\ov{\id}_H\in\Out(H)$.
In this case there is $y\in H$ such that $\ka_g=\ka_y\in\Inn(H)$,
and replacing $g\in G$ by $gy^{-1}\in G$ we still have $G=\lr{gy^{-1},H}$,
hence we may assume that $\ka_g=\id_H$, that is $g$ centralises $H$.
The cyclic extension is called {\bf faithful} if $|\ov{\ka}_g|=d$,
that is the order of $\ov{\ka}_g\in\Out(H)$ and the degree of
the extension coincide. In general, any cyclic extension 
can be written as a faithful extension of a trivial extension:

\abs
Let $|\ov{\ka}_g|=k\spmid d$. Then letting $T:=\lr{g^k,H}\unlhd G$
we have $H\leq T\leq G$ such that $G/T\cong C_k$ and 
$T/H\cong C_{\frac{d}{k}}$. Since 
$\ov{\ka}_{g^k}=(\ov{\ka}_g)^k=\ov{\id}_H\in\Out(H)$
we conclude that $T$ is a trivial extension of $H$.
In order to show that $G$ is a faithful extension of $T$,
letting $\ka_{T,g}\in\Aut(T)$ be the conjugation automorphism 
of $T$ induced by $g$, we have to verify that
$\ov{\ka}_{T,g}\in\Out(T)$ has order a multiple of $k$: 
Let $j\in\Z$ such that $(\ka_{T,g})^j\in\Inn(T)$,
thus there are $i\in\Z$ and $y\in H$ such that
$(\ka_{T,g})^j=\ka_{T,g^{ik}y}\in\Inn(T)$, implying that
$\ka_{g^{j-ik}}=\ka_y\in\Inn(H)$, hence $k\spmid j$.
\QED

\abs\abs
As is to be expected, faithful extensions are of a cohomological flavour,
and will be dealt with in the sequel \cite{ms} to the present paper.
Here, we will deal with trivial extensions, since it will turn out
that trivial GK extensions are the appropriate tool to describe 
the stems of GK trees.

\AbsT{Trivial extensions.}\label{trivextdef}
\bfa
We collect a few immediate properties of trivial extensions. 
In order to do so, let $H$ be a finite group, and let $G=\lr{g,H}$ be a 
trivial extension of $H$ of degree $d\in\N$, where $\ka_g=\id_H$.

\abs\bfi
Since $g\in G$ centralises $H$, we have $g\in Z(G)$, where the
latter denotes the centre of $G$. Moreover, we get 
$Z(H)\leq Z(G)$, or equivalently $Z(H)=Z(G)\cap H$, and thus
$Z(G)=\lr{g}Z(H)$. This yields the natural isomorphism
$$ G/Z(G)=HZ(G)/Z(G)\cong H/(Z(G)\cap H)=H/Z(H) ;$$
thus in particular we have $d=[G\cn H]=[Z(G)\cn Z(H)]$.

\abs\bfii
For the derived and lower central series
we get the following: Using the notation of \cite[Ch.1.1]{lmc},
we let $G^{(0)}:=G$ and $G^{(i+1)}=[G^{(i)},G^{(i)}]$, for $i\in\N_0$,
as well as $\gm_1(G):=G$ and $\gm_{i+1}(G)=[G,\gm_i(G)]$, for $i\in\N$,
respectively. Then, since $g\in G$ centralises $H$, 
for the derived subgroups we get $G^{(1)}=H^{(1)}$, from which we infer
$$ G^{(i)}=[G^{(i-1)},G^{(i-1)}]=[H^{(i-1)},H^{(i-1)}]=H^{(i)}, 
\quad\text{for all } i\in\N ,$$
and from $\gm_2(G)=G^{(1)}=H^{(1)}=\gm_2(H)$ we get
$$ \gm_i(G)=[G,\gm_{i-1}(G)]=[G,\gm_{i-1}(H)]=[H,\gm_{i-1}(H)]=\gm_i(H),
\quad\text{for all } i\geq 2 ;$$
thus if $H\neq\{1\}$ then $G$ has the same 
derived length and nilpotency class as $H$.

\abs\bfiii
Moreover, have the following commutative diagram,
where the horizontal arrows are the commutator maps 
associated with $G$ and $H$, respectively, and 
the vertical arrows are the natural maps induced 
by the embedding $H\leq G$, which are both isomorphisms
by the above considerations:
$$\begin{CD}
G/Z(G)\tm G/Z(G) @>{(aZ(G),bZ(G))\mt[a,b]}>> G^{(1)} \\
@A{\cong}AA @AA{\cong}A\\
H/Z(H)\tm H/Z(H) @>>{(xZ(H),yZ(H))\mt[x,y]}> H^{(1)} \\
\end{CD}$$
In particular, this implies that
$G$ and $H$ are isoclinic, see \cite[Ch.III.1]{BeylTappe}.

\AbsT{Trivial extensions of $p$-groups.}\label{trivext}
More specifically, let $G=\lr{g,H}$ be a proper trivial extension of the 
$p$-group $H$ of degree $p^l$, for some $l\in\N$, where $\ka_g=\id_H$. 
Then, by Burnside's Basis Theorem \cite[Thm.III.3.15]{hup}, 
from $G^{(1)}=H^{(1)}$ and $G^p=\lr{g^p}H^p$ we get 
$$ \Phi(G)=G^p H^{(1)}=\lr{g^p}H^p H^{(1)}=\lr{g^p}\Phi(H) .$$ 
Note that $\wti{G}:=\lr{g^p,H}=\lr{g^p}H\lhd G$, 
which is of index $p$ in $G$,
is again a trivial extension of $H$, now of degree $p^{l-1}$.
Thus, if $l\geq 2$ we get $\Phi(\wti{G})=\lr{g^{p^2}}\Phi(H)\lhd\Phi(G)$, 
which is of index $p$ in $\Phi(G)$, hence for the rank of $G$ and $\wti{G}$
we infer that $r(G)=r(\wti{G})$.

\abs
If $l=1$, that is $g^p\in H$, depending on whether $g^p\in\Phi(H)$ 
or not we have $\Phi(H)=\Phi(G)$ or $\Phi(H)<\Phi(G)$, respectively.
Thus in the former case we get $r(G)=r(H)+1$,
while in the latter case, since $g^{p^2}\in\Phi(H)$ anyway, we infer
that $\Phi(H)$ has index $p$ in $\Phi(G)$, hence $r(G)=r(H)$.

\abs
In particular, if $G$ is a trivial GK extension of its root $R$
of level $l\in\N$, then 
$$ 1\leq r(G)=r(\cK(G))=\cdots=r(\cK^{l-1}(G))\leq r(R)+1 .$$ 

\AbsT{Parents.}\label{parents}
We are now prepared to put things into a more structural perspective,
showing that the trivial extensions of a fixed finite group $H$ 
of degree $d$ coincide with certain epimorphic images 
of suitable universal groups:

\abs\bfa
To this end, with some foresight we fix $h\in Z(H)$, where
the case $h=1$ is explicitly allowed.
Moreover, let $C_{d\cdot |h|}$ 
be an abstract copy of the cyclic group of order $d\cdot |h|$, 
and let $\wh{g}\in C_{d\cdot |h|}$ be a generator, that is we have
$\lr{\wh{g}}=C_{d\cdot |h|}$. Then the associated {\bf parent group} 
is defined as the direct product
$$ \wh{G}:=\lr{\wh{g}}\tm H .$$
Now we consider the centrally amalgamated product $\wh{G}/C$, where
$$ C:=\lr{(\wh{g}^d,h^{-1})}\leq\lr{\wh{g}}\tm Z(H)=Z(\wh{G})\leq\wh{G} ,$$ 
a cyclic central subgroup of order $|h|$.
Note that $\wh{G}$ only depends on the order $|h|$,
while $\wh{G}/C$ explicitly depends on the choice of $h$; 
in particular we have $C=\{1\}$ if and only if $h=1$.
We show that $\wh{G}/C$ is a trivial extension of $H$ 
of degree $d$ such that $\ka_{\wh{g}}=\id_H$:

\abs
By construction, we have $(\lr{\wh{g}}\tm\{1\})\cap C=\{1\}$
and $(\{1\}\tm H)\cap C=\{1\}$, hence both $\lr{\wh{g}}$ and $H$ 
can be considered as normal subgroups of $\wh{G}/C$, 
via the natural embeddings
$\lr{\wh{g}}\ra\wh{G}/C\cn\wh{g}\mt(\wh{g},1)C$ and
$H\ra\wh{G}/C\cn y\mt (1,y)C$, respectively.
Using these identifications, we have 
$\wh{G}/C=\lr{\wh{g},H}$, where $\wh{g}$ centralises $H$.
Moreover, for $k\in\Z$ and $y\in H$ we have 
$(\wh{g},1)^kC=(1,y)C\in\wh{G}/C$
if and only if $(\wh{g}^k,y^{-1})\in C\leq\wh{G}$,
which holds if and only if $d\spmid k$ and $y=h^{\frac{k}{d}}$.  
Hence, again using the above identifications, we get
$\lr{\wh{g}}\cap H=\lr{\wh{g}^d}=\lr{h}\leq\wh{G}/C$,
implying the assertion concerning the degree;
note that this also shows that $\wh{g}^d=h\in H\leq\wh{G}/C$.
\QED

\abs\bfb
This shows that the quotients of $\wh{G}$ with respect to
the cyclic central subgroups of the form considered above all
are trivial extensions of $H$ of degree $d$;
in particular, the case $h=1$ shows that $\wh{G}$ itself
is such an extension of $H$. We finally show that 
any trivial extension of $H$ of degree $d$ can be realised like this:

\abs
Let $G=\lr{g,H}$ be a trivial extension of $H$ 
of degree $d$, where $\ka_g=\id_H$. Then let
$h:=g^d\in Z(H)$, hence we have $\lr{g}\cap H=\lr{h}$
and thus $|g|=d\cdot |h|$. Now consider the parent group 
$\wh{G}:=\lr{\wh{g}}\tm H$, where $|\wh{g}|=|g|$.
Hence, since $g\in G$ centralises $H$, there is an epimorphism
$\wh{G}\ra G\cn(\wh{g}^k,y)\mt g^k\cdot y$, for all $k\in\Z$ and $y\in H$.
Since we have $g^k\cdot y=1$, or equivalently 
$g^k=y^{-1}\in\lr{g}\cap H=\lr{h}$, 
if and only if $d\spmid k$ and $y=h^{-\frac{k}{d}}$,
we infer that the above epimorphism has cyclic central kernel
$C:=\lr{(\wh{g}^d,h^{-1})}\unlhd\wh{G}$,
implying that $G\cong\wh{G}/C$ as desired.
\QED

\AbsT{Automorphisms of parents.}\label{automdirprod}
We collect a few facts on automorphism
groups of parent groups, where in view of the applications
to come we impose a further technical restriction:
More precisely, let still $H$ be a finite group, 
let $C_n$ be an abstract copy of the cyclic group of order $n\in\N$,
where now
$$ \exp(Z(H))\spmid n ,$$
let $\wh{g}$ be a generator, that is we have $\lr{\wh{g}}=C_n$,
and let $\wh{G}:=\lr{\wh{g}}\tm H$ be the associated parent group.
We consider the group of automorphisms 
$$ \Aut_0(\wh{G}):=\{\wh{\ph}\in\Aut(\wh{G});H^{\wh{\ph}}=H\}
   \leq\Aut(\wh{G}) :$$

\abs
Recall that both $\lr{\wh{g}}$ and $H$ are identified 
with subgroups of $\wh{G}$.
Since $H$ is invariant under any $\wh{\ph}\in\Aut_0(\wh{G})$,
we have a restriction homomorphism
$\Aut_0(\wh{G})\ra\Aut(H)\cn\wh{\ph}\mt\wh{\ph}|_H$,
and the natural epimorphism $\wh{G}\ra\wh{G}/H$ induces 
a homomorphism 
$\Aut_0(\wh{G})\ra\Aut(\wh{G}/H)\cong\Aut(C_n)\cong\Z_n^\ast$.

\abs
Hence, given $\wh{\ph}\in\Aut_0(\wh{G})$, let 
$\al\in\Aut(H)$ be its restriction to $H$,
and let $k\in\Z_n^\ast$ and $z\in H$ such that 
$\wh{g}^{\wh{\ph}}=(\wh{g}^k,z)$,
where since $\wh{g}\in Z(\wh{G})$ we infer that 
$\wh{g}^{\wh{\ph}-k}:=\wh{g}^{\wh{\ph}}\cdot\wh{g}^{-k}
=z\in Z(\wh{G})\cap H=Z(H)$.
Thus $\wh{\ph}$ is uniquely described by 
$$ (\al,k,z)\in\Aut(H)\tm\Z_n^\ast\tm Z(H) .$$

\abs
Conversely, given such a triple $(\al,k,z)\in\Aut(H)\tm\Z_n^\ast\tm Z(H)$,
since we have $|z|\spmid\exp(Z(H))\spmid n=|\wh{g}|$ there is the
homomorphism $\lr{\wh{g}}\ra Z(H)\cn \wh{g}\mt z$, 
and thus we get an automorphism $\wh{\ph}\in\Aut_0(\wh{G})$ by letting
$$ (\wh{g}^i,y)^{\wh{\ph}}:=(\wh{g}^{i\cdot k},z^i\cdot y^\al),
\quad\text{for all }i\in\Z,\, y\in H .$$ 

\abs
Given $\wh{\ph}'\in\Aut_0(\wh{G})$ with associated triple 
$(\al',k',z')\in\Aut(H)\tm\Z_n^\ast\tm Z(H)$, we have
$$ \wh{g}^{\wh{\ph}\wh{\ph}'}=(\wh{g}^k,z)^{\wh{\ph}'}
 =(\wh{g}^{kk'},(z')^k\cdot z^{\al'})
\quad\text{and}\quad
y^{\wh{\ph}\wh{\ph}'}=y^{\al\al'} ,$$ 
for all $y\in H$.
Hence $\wh{\ph}\wh{\ph}'$ is associated with the triple
$(\al\al',kk',(z')^k\cdot z^{\al'})\in\Aut(H)\tm\Z_n^\ast\tm Z(H)$.
Thus in conclusion we have  
$$ \Aut_0(\wh{G})\cong(\Aut(H)\tm\Z_n^\ast)\ltimes Z(H) ,$$
where $\Aut(H)\tm\Z_n^\ast$ acts on $Z(H)$ by
$$ Z(H)\ra Z(H)\cn z\mt (z^\al)^{k^{-1}}=(z^{k^{-1}})^\al,
\quad\text{for all }\al\in\Aut(H),\, k\in\Z_n^\ast .$$

%% file: stems.tex
\section{Stems}\label{stems}

\abs
We are now prepared to state and prove a group theoretic criterion
saying whether a group not of GK type has an infinite tree attached
to it, and in this case describing the branching behaviour of
its stems and how the groups on the stems look like.
As is to expected, trivial extension will play the crucial role.

\Thm\label{trivextthm}
Let $H$ be a finite $p$-group. Then $H$ has a proper trivial 
GK extension if and only if $\exp(Z(H))=\exp(H)$. 

\abs
In this case, $H$ has proper trivial GK extensions of any $p$-power degree, 
and for any such extension $G$ we have $\exp(Z(G))=\exp(G)$.

\Pf
Let $G$ be a proper trivial GK extension of $H$, and let 
$g\in G$ such that $G=\lr{g,H}$ and $\ka_g=\id_H$. Hence we have
$G/H\cong C_{p^l}$ for some $l\in\N$, and letting $\exp(H)=p^e$,
for some $e\in\N_0$, we have $\exp(G)=p^{e+l}$. Then since $H\leq\cK(G)$
we have $g\in G\smin\cK(G)$, implying that $|g|=p^{e+l}$.
Hence, since $g\in Z(G)$ we have $\exp(Z(G))=\exp(G)$.
Moreover, $h:=g^{p^l}\in Z(G)\cap H=Z(H)$, see \ref{trivextdef},
has order $|h|=p^e$, and hence we have $\exp(Z(H))=\exp(H)$. 

\abs
Let conversely $\exp(Z(H))=\exp(H)=p^e$, for some $e\in\N_0$, and let
$h\in Z(H)$ such that $|h|=p^e$. 
For any $l\in\N$ let $\wh{G}:=\lr{\wh{g}}\tm H$ be the parent group
where $|\wh{g}|=p^{e+l}$, and let $\wh{G}/C$ be the 
centrally amalgamated product with respect to
$C:=\lr{(\wh{g}^{p^l},h^{-1})}\unlhd\wh{G}$.
By \ref{parents}, the group $\wh{G}/C$ is a trivial extension
of $H$ of degree $p^l$, where $\wh{g}^{p^l}=h\in\wh{G}/C$
and $\ka_{\wh{g}}=\id_H$.

\abs
Since $\exp(H)=p^e\spmid p^{e+l}=|\wh{g}|$ and $\wh{g}\in Z(\wh{G}/C)$, 
we infer that $\exp(\wh{G}/C)=|\wh{g}|=p^{e+l}$.
Moreover, for all $y\in H$ and $i\in\{0\ld l-1\}$,
if $k\in\Z$ such that $p^i\spmid k$, but $p^{i+1}\spnmid k$,
then we have $|\wh{g}^k y|=p^{e+l-i}$,
while if $p^l\spmid k$ then we get $|\wh{g}^k y|\spmid p^e$.
Hence for $i\in\{1\ld l\}$ we have
$$ \cK^i(\wh{G}/C):=\{y\in\wh{G}/C;|y|\spmid p^{e+l-i}\}
 =\lr{\wh{g}^{p^i},H}\unlhd\wh{G}/C ,$$
in particular we get $\cK^l(\wh{G}/C)=H$.
Letting $\cK^0(\wh{G}/C):=\wh{G}/C$, this yields
$$ \cK^{i-1}(\wh{G}/C)/\cK^i(\wh{G}/C)\cong
   \lr{\wh{g}^{p^{i-1}}}/\lr{\wh{g}^{p^{i}}}\cong C_p ,$$
for all $i\in\{1\ld l\}$, implying that $\wh{G}/C$ is of GK type,
where since $H$ occurs in the GK series of $\wh{G}/C$,
the latter is a GK extension of $H$.
\QED

\Cor\label{trivextcor}
Let $(p^{e_1}\ld p^{e_r})$, where $r\in\N_0$ and 
$0<e_1\leq\cdots\leq e_r=e$, be the abelian invariants of $Z(H)$;
we allow for $r=0$, letting $e_0:=0$, to catch the case $H=Z(H)=\{1\}$.
Moreover, let $G=\lr{g,H}$ be a proper trivial 
GK extension of $H$ of degree $p^l$, where $l\in\N$, such that $\ka_g=\id_H$.

\abs
Then $Z(G)=\lr{g}Z(H)$ is a proper trivial GK extension of $Z(H)$ 
of degree $p^l$, where $\cK_i(Z(G))=\lr{g^{p^i}}Z(H)$ for $i\in\{1\ld l\}$,
and its abelian invariants are 
$$ (p^{e_1}\ld p^{e_{r-1}},p^{e_r+l}) .$$

\Pf
We only have to show the assertion on abelian invariants:
Writing $Z(H)=\prod_{i=1}^r\lr{h_i}$, where $|h_i|=p^{e_i}$, 
since $p^{e_r}=p^e=\exp(Z(H))$
we may choose $h_r:=g^{p^l}\in Z(H)$, implying
$Z(G)=\left(\prod_{i=1}^{r-1}\lr{h_i}\right)\tm\lr{g}$.
\QED

\Thm\label{inftree}
Let $R$ be a finite $p$-group not of GK type. 

\abs\bfa
Then the following statements are equivalent:

\bfi
The GK tree $\cT(R)$ is infinite, or equivalently $\cT(R)$ has a stem.

\bfii
The group $R$ has a (proper) trivial GK extension.
 
\bfiii
We have $\exp(R)=\exp(Z(R))$.

\abs\bfb
If $\cT(R)$ is infinite, then for any GK extension $G$ of $R$ 
the following statements are equivalent:

\bfi
The group $G$ lies on a stem of $\cT(R)$.

\bfii
The group $G$ is a trivial extension of $R$.

\bfiii
We have $\exp(G)=\exp(Z(G))$.

\Pf
\bfa
The equivalence of (ii) and (iii) has been shown in \ref{trivextthm};
note that since $R$ is not of GK type, properness is automatic.
If $\exp(R)=\exp(Z(R))$, then again by \ref{trivextthm} $R$ has trivial 
GK extensions of arbitrarily large degree, hence
$\cT(R)$ is infinite, showing that (iii) implies (i).

\abs
Let finally $\cT(R)$ be infinite, and let $H$ be a GK extension of $R$ 
of level $k\in\N$ lying on a stem of $\cT(R)$.
Hence letting $p^f=\exp_p(\Aut(R))$ be the $p$-exponent of $\Aut(R)$,
where $f\in\N_0$, there is a GK extension $G$ of $R$ of level $l:=k+f$
such that $\cK^f(G)=H$. Thus we have
$$ R=\cK^k(H)=\cK^l(G)<H=\cK^f(G)\leq G .$$
Then for any $g\in G\smin\cK(G)$
we have $g^{p^f}\in\cK^f(G)\smin\cK^{f+1}(G)=H\smin\cK(H)$ and 
$(\ka_g)^{p^f}=\id_R\in\Aut(R)$, hence we infer that $H$ is a 
trivial GK extension of $R$, showing that (i) implies (ii).

\abs\bfb
The above argument shows that any GK extension $G$ of $R$ 
lying on a stem of $\cT(R)$ is a trivial extension,
that is (i) implies (ii). Moreover, if $G$ is a trivial 
extension of $R$, then by \ref{trivextthm} 
we have $\exp(G)=\exp(Z(G))$, showing that (ii) implies (iii).
Finally, if $\exp(G)=\exp(Z(G))$ then 
once again by \ref{trivextthm} the group $G$ has trivial 
GK extensions of arbitrarily large $p$-power degree, 
implying that $G$ lies on a stem of $\cT(R)$,
showing hat (iii) implies (i).
\QED

\AbsT{Action on centres.}\label{actcen}
Let $H$ be a finite $p$-group such that $\exp(Z(H))=\exp(H)=p^e$ 
for some $e\in\N_0$.
We proceed to describe the isomorphism types of trivial GK extensions
of $H$. To this end, we need some preparations first:

\abs
For $i\in\N_0$ let $Z^{p^i}(H):=(Z(H))^{p^i}\unlhd Z(H)$, where
$Z(H)$ being abelian implies
$$ Z^{p^i}(H)=\{z^{p^i}\in Z(H);z\in Z(H)\} .$$
Hence we have $Z^1(H)=Z(H)$ and $Z^{p^{i+1}}(H)=(Z^{p^i}(H))^{p}$,
for $i\in\N_0$, and for $i\in\{0\ld e\}$
we get a strictly descending chain of characteristic subgroups
$$ \{1\}=Z^{p^e}(H)<Z^{p^{e-1}}(H)<\cdots<Z^p(H)<Z^1(H)=Z(H) .$$
Moreover, let 
$$ \cZ(H):=\{z\in Z(H);|z|=p^e\} ,$$
where for $e>0$, that is $H\neq\{1\}$, we have $\cZ(H)\sseq Z(H)\smin Z^p(H)$,
while of course $Z(\{1\})=Z^p(\{1\})=\cZ(\{1\})$.
Hence the condition in \ref{trivextthm} is equivalent to saying
$\cZ(H)\neq\emp$.
Now there are various groups acting on $\cZ(H)$:

\abs\bfi
Since for all $z\in\cZ(H)$ and $y\in Z^{p^i}(H)$, where $i\geq 1$,
we have $|zy|=p^e=|z|$, we conclude that $Z^{p^i}(H)$ acts faithfully,
even semi-regularly, on $\cZ(H)$ by translation
$$ 
   \cZ(H)\ra\cZ(H)\cn z\mt zy,
   \quad\text{where }y\in Z^{p^i}(H) .$$
\bfii
Moreover, $\Z_{p^e}^\ast$ acts faithfully,
even semi-regularly, on $\cZ(H)$ by exponentiation
$$ 
   \cZ(H)\ra\cZ(H)\cn z\mt z^k,
   \quad\text{where }k\in\Z_{p^e}^\ast .$$
\bfiii
Finally, $\Aut(H)$ acts on $\cZ(H)\sseq Z(H)$ by the natural action
$$ 
   \cZ(H)\ra\cZ(H)\cn z\mt z^\al,
   \quad\text{where }\al\in\Aut(H) ;$$
note that, since $\Inn(H)$ acts trivially on $Z(H)$, this action
factors through $\Out(H)$, but even the action of $\Out(H)$ 
need not be faithful.

\abs
Thus, in combination, for $i\geq 1$
the action of $(\al,k,y)\in\Aut(H)\tm\Z_{p^e}^\ast\tm Z^{p^i}(H)$ 
on $\cZ(H)$ is given as
$$ \cZ(H)\ra\cZ(H)\cn z\mt (z^k)^{\al}\cdot y=(z^{\al})^k\cdot y ,$$
and concatenation with 
$(\al',k',y')\in\Aut(H)\tm\Z_{p^e}^\ast\tm Z^{p^i}(H)$ yields 
$$ \cZ(H)\ra\cZ(H)\cn z\mt (z^{kk'})^{\al\al'}\cdot(y^{k'})^{\al'}\cdot y' 
                          =(z^{kk'})^{\al\al'}\cdot(y^{\al'})^{k'}\cdot y'.$$
Hence we have an action homomorphism 
$$ (\Aut(H)\tm\Z_{p^e}^\ast)\ltimes Z^{p^i}(H)\ra\cS_{\cZ(H)} ,$$
where again the action of $\Aut(H)$ on $Z^{p^i}(H)$ is induced by the
natural action of $\Aut(H)$ on $Z^{p^i}(H)\sseq Z(H)$.
We let $A_i(H)\leq\cS_{\cZ(H)}$ be its image, that is 
the permutation group on $\cZ(H)$ generated by this action. 
For completeness, let $A_0(H):=\cS_{\cZ(H)}$ just be the full
symmetric group on $\cZ(H)$.

\abs
Since for $i\geq e\geq 1$ we have $Z^{p^i}(H)=\{1\}$,
we infer that $A_i(H)$ is an epimorphic image of
$\Aut(H)\tm\Z_{p^e}^\ast$, yielding a chain of normal subgroups
$$ \cdots=A_{e+1}(H)=A_e(H)\leq A_{e-1}(H)\leq\cdots\leq A_1(H) ;$$
note that for $e=0$, that is $H=\{1\}$, we get 
$A_i(H)=\{1\}$ for all $i\in\N_0$ anyway.

\abs
Finally, in view of the subsequent result
in conjunction with the comments on the rank of trivial extensions
in \ref{trivext}, we note that 
for any $A_l(H)$-orbit $\cO\sseq\cZ(H)$, where $l\in\N$,
we indeed have either $\cO\sseq\Phi(H)$ or $\cO\cap\Phi(H)=\emp$:

\abs
It is immediate that $\cZ(H)\cap\Phi(H)$ is invariant under
both the exponentiation action of $\Z_{p^e}^\ast$, and
the natural action of $\Aut(H)$, and since $Z^p(H)\leq H^p\leq\Phi(H)$, 
by Burnside's Basis Theorem \cite[Thm.III.3.15]{hup},
the same holds with respect to the translation action of $Z^{p^l}(H)$.
\QED

\Thm\label{stemthm}
Let $H$ be a finite $p$-group such that $\exp(Z(H))=\exp(H)=p^e$
for some $e\in\N_0$.
Then the (isomorphism types of) proper trivial GK extensions of $H$ 
of degree $p^l$, where $l\in\N$, are in bijection with the
$A_l(H)$-orbits in $\cZ(H)$.

\abs
More precisely, the bijection is given by mapping $h\in\cZ(H)$ 
to $\wh{G}/C$, where $\wh{G}=\lr{\wh{g}}\tm H$ is the parent
group such that $|\wh{g}|=p^{e+l}$, and 
$C:=\lr{(\wh{g}^{p^l},h^{-1})}\unlhd\wh{G}$. 

\Pf
By the proof of \ref{trivextthm}, the proper trivial GK extensions of $H$
are quotients of the parent group $\wh{G}$ as specified above. Hence
let $\wh{G}/C$ and $\wh{G}/C'$ be trivial GK extensions of $H$ with 
respect to $h\in\cZ(H)$ and $h'\in\cZ(H)$, respectively, that is we have
$$ C:=\lr{(\wh{g}^{p^l},h^{-1})}\unlhd\wh{G}
\quad\text{and}\quad
   C':=\lr{(\wh{g}^{p^l},(h')^{-1})}\unlhd\wh{G} ,$$
and let $\psi\cn\wh{G}/C\ra\wh{G}/C'$ be an isomorphism.
We have to show that $h$ and $h'$ are in the same $A_l(H)$-orbit:

\abs
From $HC/C=\{x\in\wh{G}/C;|x|\spmid p^e\}$ and 
$HC'/C'=\{x\in\wh{G}/C';|x|\spmid p^e\}$
we conclude that $\psi$ restricts to an isomorphism 
$H\cong HC/C\ra HC'/C'\cong H$, 
which by the usual identifications yields an automorphism $\al\in\Aut(H)$.
Moreover, letting $g:=\wh{g}C\in\wh{G}/C$ and $g':=\wh{g}C'\in\wh{G}/C'$, 
we have $g^\psi=(g')^k\cdot z\in\wh{g}C'$, for some $k\in\Z$ and $z\in H$.
Since $|(g')^k\cdot z|=|g|=|\wh{g}|=p^{e+l}$ and $\wh{g}\in Z(\wh{G})$,
we conclude that $k\in\Z_{p^{e+l}}^\ast$ and $z\in Z(H)$. 
Moreover, from $h=g^{p^l}\in\wh{G}/C$ and $h'=(g')^{p^l}\in\wh{G}/C'$  
we get
$$ h^\al=(g^{p^l})^\psi=(g')^{p^l\cdot k}\cdot z^{p^l}
 =(h')^k\cdot z^{p^l} ,$$
or equivalently
$$ h'=(h^\al\cdot z^{-p^l})^{k^{-1}}
     =(h^\al)^{k^{-1}}\cdot (z^{-k^{-1}})^{p^l}\in Z(H) .$$
Hence $h\in\cZ(H)$ is mapped to $h'\in\cZ(H)$ by the action of
$(\al,\ov{k}^{-1},(z^{-\ov{k}^{-1}})^{p^l})
 \in\Aut(H)\tm\Z_{p^e}^\ast\tm Z^{p^l}(H)$,
hence $h$ and $h'$ indeed are in the same $A_l(H)$-orbit; 
note that the exponentiation action of $\Z_{p^{e+l}}^\ast$ on $Z(H)$
factors through the action of $\Z_{p^e}^\ast$ 
via the natural epimorphism $\baraut\cn\Z_{p^{e+l}}^\ast\ra\Z_{p^e}^\ast$.

\abs
Conversely, let $h\in\cZ(H)$ and $C:=\lr{(\wh{g}^{p^l},h^{-1})}\unlhd\wh{G}$,
and for some 
$(\al,\ov{k},z^{p^l})\in\Aut(H)\tm\Z_{p^e}^\ast\tm Z^{p^l}(H)$,
let $h':=(h^\al)^k\cdot z^{p^l}\in\cZ(H)$ and 
$C':=\lr{(\wh{g}^{p^l},(h')^{-1})}\unlhd\wh{G}$.
We have to show that $\wh{G}/C$ and $\wh{G}/C'$ are isomorphic: 

\abs
To this end, by \ref{automdirprod} let $\wh{\psi}\in\Aut_0(\wh{G})$ be the
automorphism of $\wh{G}$ described by the triple 
$(\al,z^{-k^{-1}},k^{-1})\in\Aut(H)\tm Z(H)\tm\Z_{p^{e+l}}^\ast$.
Then in $\wh{G}$ we have
$$ ((\wh{g}^{p^l},h^{-1})^{\wh{\psi}})^k
  =(\wh{g}^{p^l\cdot k^{-1}},z^{-p^l\cdot k^{-1}}\cdot h^{-\al})^k
  =(\wh{g}^{p^l},z^{-p^l}\cdot (h^{\al})^{-k})
  =(\wh{g}^{p^l},(h')^{-1}) ,$$
showing that $C^{\wh{\psi}}=C'$, and thus $\wh{\psi}$ induces an
isomorphism $\wh{G}/C\ra\wh{G}/C'$.
\QED

\AbsT{Stems of GK trees.}\label{stemrem}
In conclusion, we are now able to describe 
the stems of infinite GK trees and their branching behaviour:
Let $R$ be a finite $p$-group not of GK type, being the root
of an infinite GK tree, that is we have $\exp(Z(R))=\exp(R)=p^e$, 
for some $e\in\N_0$. 

\abs
If $H$ lies on a stem of $\cT(R)$ and has level $l\in\N_0$,
then the number of stems into which this stem branches at $H$,
that is the number of groups $G$ lying on stems of 
$\cT(R)$ such that $\cK(G)=H$, coincides with 
the number of $A_{l+1}(R)$-orbits into which the
$A_l(R)$-orbit in $\cZ(R)$ associated with $H$, in the sense of
\ref{stemthm}, splits; note that since we have agreed on letting
$A_0=\cS_{\cZ(R)}$, which is transitive on $\cZ(R)$,
this in particular holds for the level $l=0$.

\abs
Hence, since $A_{l+1}(R)=A_l(R)$ as soon as $l\geq e$,
branching occurs at most at levels $l\in\{0\ld e-1\}$.
In particular, for the trivial group $R=\{1\}$, that is $e=0$,
we recover the result that $\cT(\{1\})$ has a single stem; see \ref{trivgrp}.
Moreover, in general, the total number of stems of $\cT(R)$ 
coincides with the number of $A_e(R)$-orbits in $\cZ(R)$,
thus is finite.

\abs 
In view of the examples given in \ref{stem2ex} and \ref{branchexmpl}, 
it seems that stronger general statements concerning the branching 
behaviour of stems of infinite GK trees cannot be hoped for.
Actually, we are better off for our favourite example
of abelian root groups:

\ExT{Abelian groups}\label{abelianII}
Let $R$ be an abelian $p$-group not of GK type.
Then by \ref{inftree} the associated GK tree $\cT(R)$ is infinite.
Moreover, the groups lying on a stem of $\cT(R)$ being
precisely the trivial GK extensions of $R$,
we infer from \ref{trivextcor} that a group $G$ in $\cT(R)$
lies on a stem if and only if $G$ is abelian.
In that case, the abelian invariants of $G$ are uniquely determined
by the GK level of $G$, implying that $\cT(R)$ has a unique stem.

\abs
The latter statement can, alternatively, also be seen as follows:
Since $R$ is abelian, any cyclic subgroup of $R$ of maximal order
has a complement in $R$, and any bijective exponentiation map is an
automorphism of $R$. This implies that $\Aut(R)$ acts
transitively on the set $\cZ(R)=\{z\in R;|z|=\exp(R)\}$.
Thus we conclude that $\cZ(R)$ consists of a single $A_l(R)$-orbit,
for all $l\in\N_0$, hence the uniqueness of the stem of $\cT(R)$
also follows from \ref{stemrem}.

%% file: 2grpex.tex
\section{Examples of even order}\label{2grpex}

\abs
We conclude the paper with an extended collection of explicit examples.
Although we try to give specific theoretical descriptions, 
we point out that the examples typically have been found initially
by searching the {\sf SmallGroups} database \cite{SmGrp},
which is available through the computer algebra system \GAP{} \cite{GAP}, 
and contains all the finite groups up to order $1023$ (and many more).

\abs
Moreover, using the facilities to compute with $p$-groups 
available in \GAP{}, it is straightforward to check whether
a given group is of GK type, in this case to determine its
GK series, as well as the associated subset $\cZ$ of its centre, 
and the orbits of the action of the permutation groups $A_l$ on it.

\AbsT{The trees rooted at small $2$-groups.}\label{2grpexmpl}
The GK trees rooted in $2$-groups of order at most $16$, and 
having more than one vertex, are given in Tables 
\ref{treetbl1}--\ref{treetbl3}.
The stems of the trees are drawn vertically, while bushes
branch off to the right.
Left to the trees we indicate the order of the groups
in the various layers, thus edges are directed downwards.
Attached to each vertex we give the number of the associated
group in the {\sf SmallGroups} database.
Moreover, at the bottom of the trees we indicate the rank
of the groups in the various columns, in order to illustrate
the statements in \ref{immediate} and \ref{trivext}.
Finally, for the trees in Table \ref{treetbl1} and Table \ref{treetbl2},
which are rooted at abelian groups, we also indicate the 
isomorphism types of the groups lying on the stem; recall that 
by \ref{abelianII} any tree with an abelian root has a unique stem.

\abs
Usually, we cannot resist to draw the trees a bit further than
is justified by the existing data, in order to
indicate their infiniteness and to point out the expected
periodic behaviour. In the tables, the trees are proven 
to be correct for all layers carrying a group order,
and the existence of stems can be derived from \ref{stemthm}. 
But since so far we do not have a theory describing the bushes, 
which is left to the sequel \cite{ms} of the present paper,
apart from that they are (mostly) conjectural,
with the exception of the trees $\cT(\{1\})$ and $\cT(C_p^2)$ for
arbitrary prime $p$, given in Table \ref{treetbl1} and Table 
\ref{treegentbl1}, whose correctness is proved in \ref{trivgrp} 
and \ref{cp2tree}, respectively.

\abs
More precisely, for the non-abelian groups of order $8$ we get the following:
For the dihedral group $D_8\cong\textsf{SmallGroups}(8,3)$
the set $\cK(D_8)=\{x\in D_8;|x|\leq 2\}$ has precisely six elements, 
thus is not a subgroup;
for the quaternion group $Q_8\cong\textsf{SmallGroups}(8,4)$
we have $\cK(Q_8)=\{x\in Q_8;|x|\leq 2\}=Z(Q_8)$,
which is a non-maximal subgroup.
Hence both groups are not of GK type, and it turns out that
neither of them occurs as a root of a group of GK type. 

\abs
From the $14$ (isomorphism types of)
groups of order $16$ there are five of GK type,
and thus occur in Table \ref{treetbl1};
another five are not of GK type but are roots of groups of GK type,
and their trees are given in Table \ref{treetbl2} and Table \ref{treetbl3},
for the abelian and non-abelian cases, respectively; 
and the remaining four are neither of GK type nor roots of groups of GK type.
We now have a closer look at the groups rooting the trees in 
Table \ref{treetbl3}: 

\begin{table}\caption{Trees rooted at groups of order dividing $8$.}
\label{treetbl1}

\abs
\includegraphics[height=80mm]{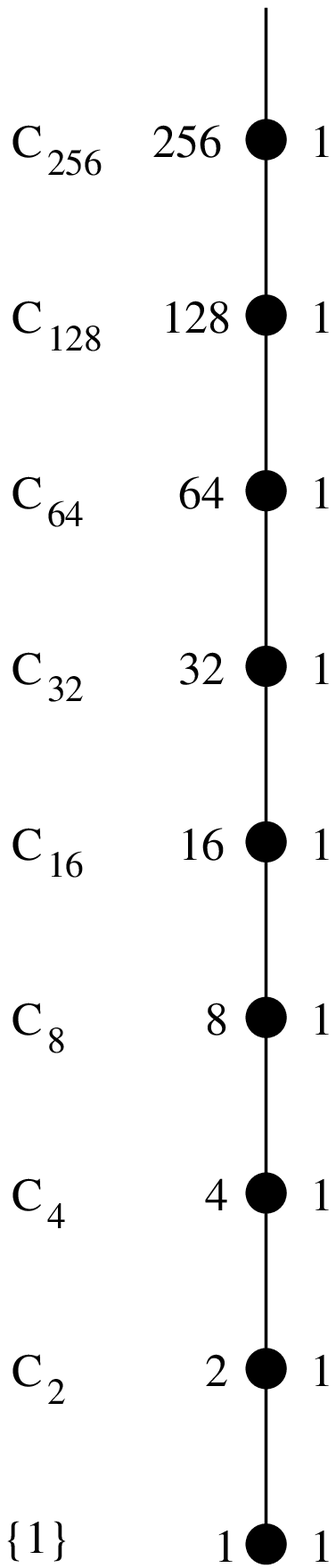} 
\hfill
\includegraphics[height=80mm]{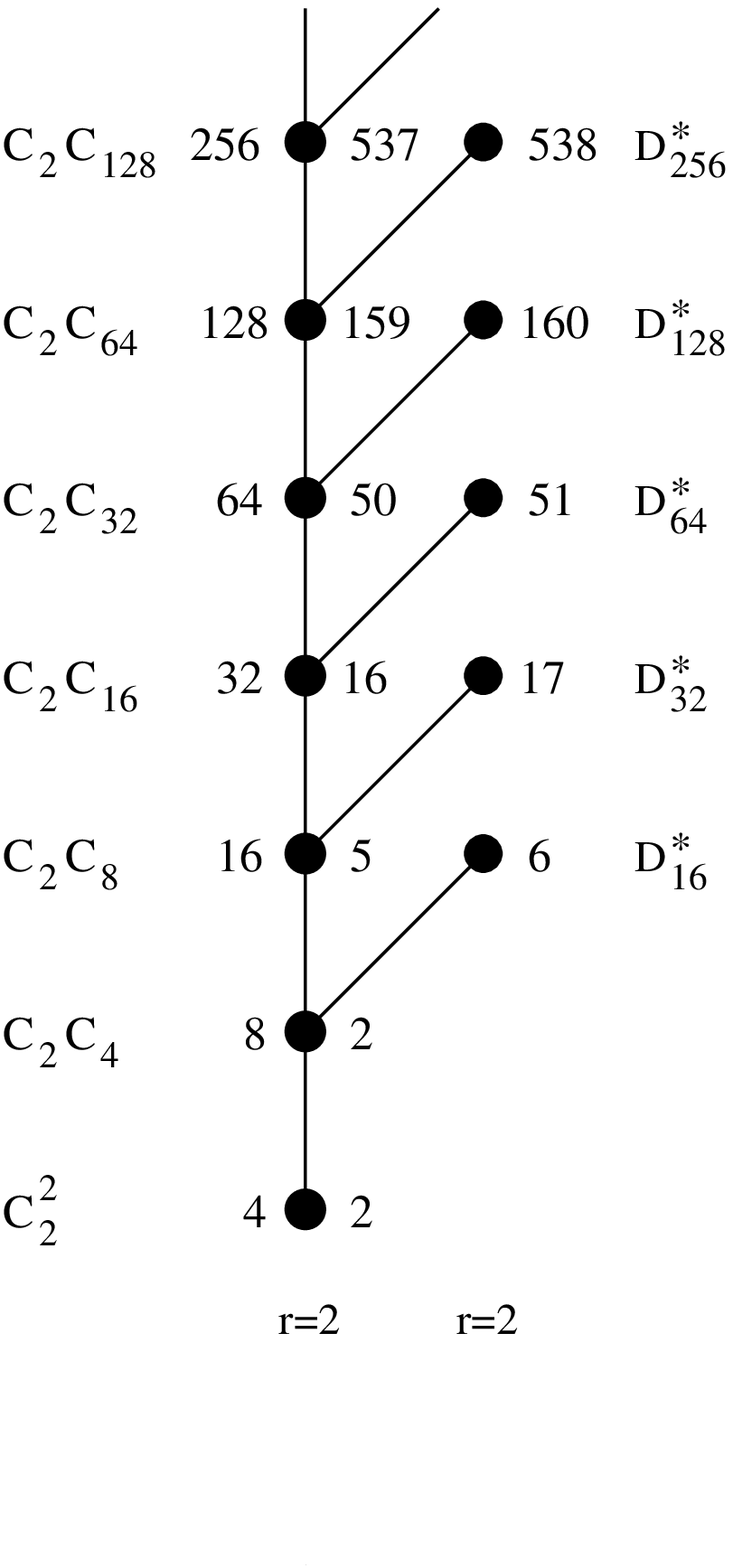} 
\hfill
\includegraphics[height=80mm]{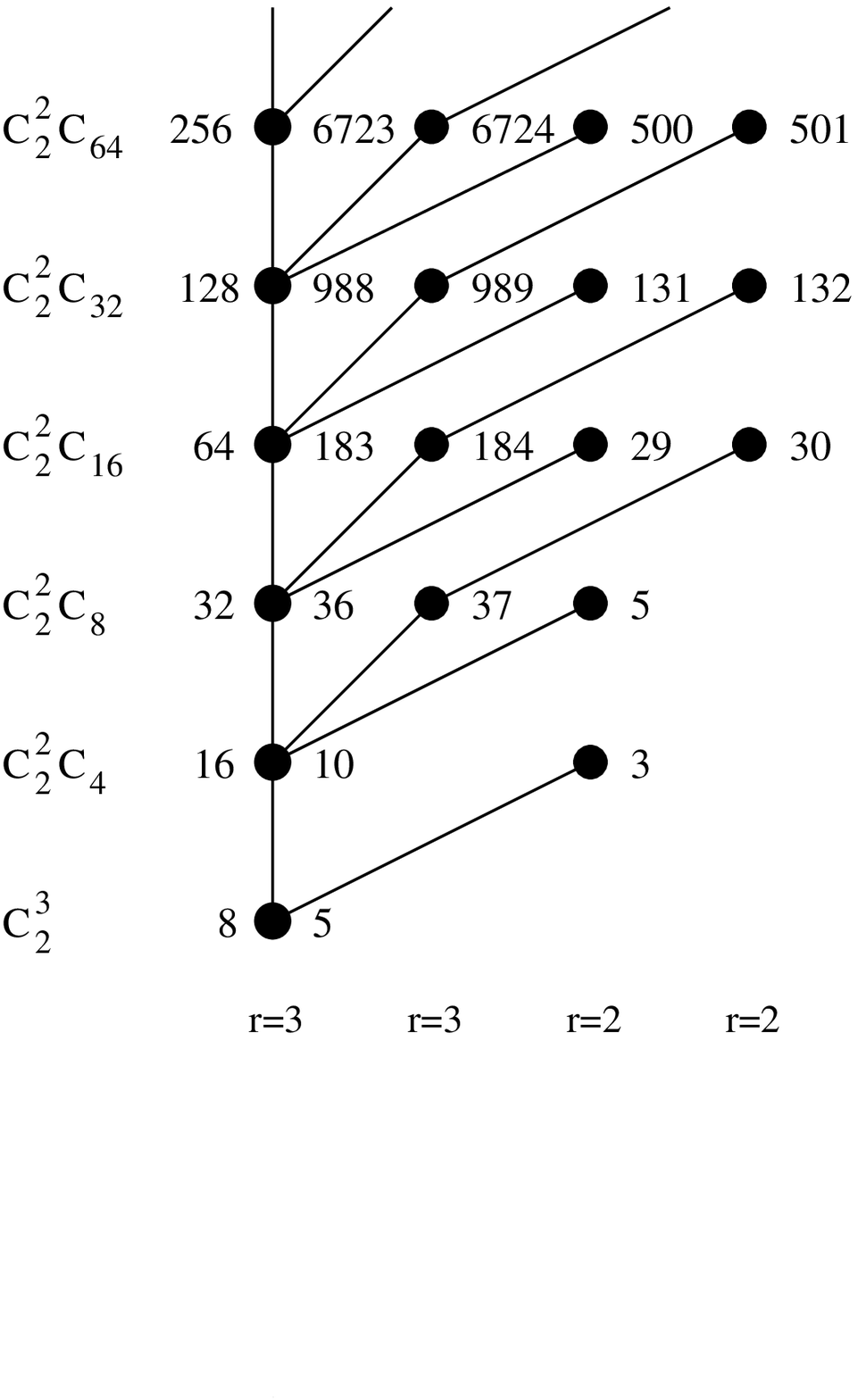} 

\abs\hrulefill
\end{table}

\begin{table}\caption{Trees rooted at abelian groups of order $16$.}
\label{treetbl2}

\abs
\includegraphics[height=42mm]{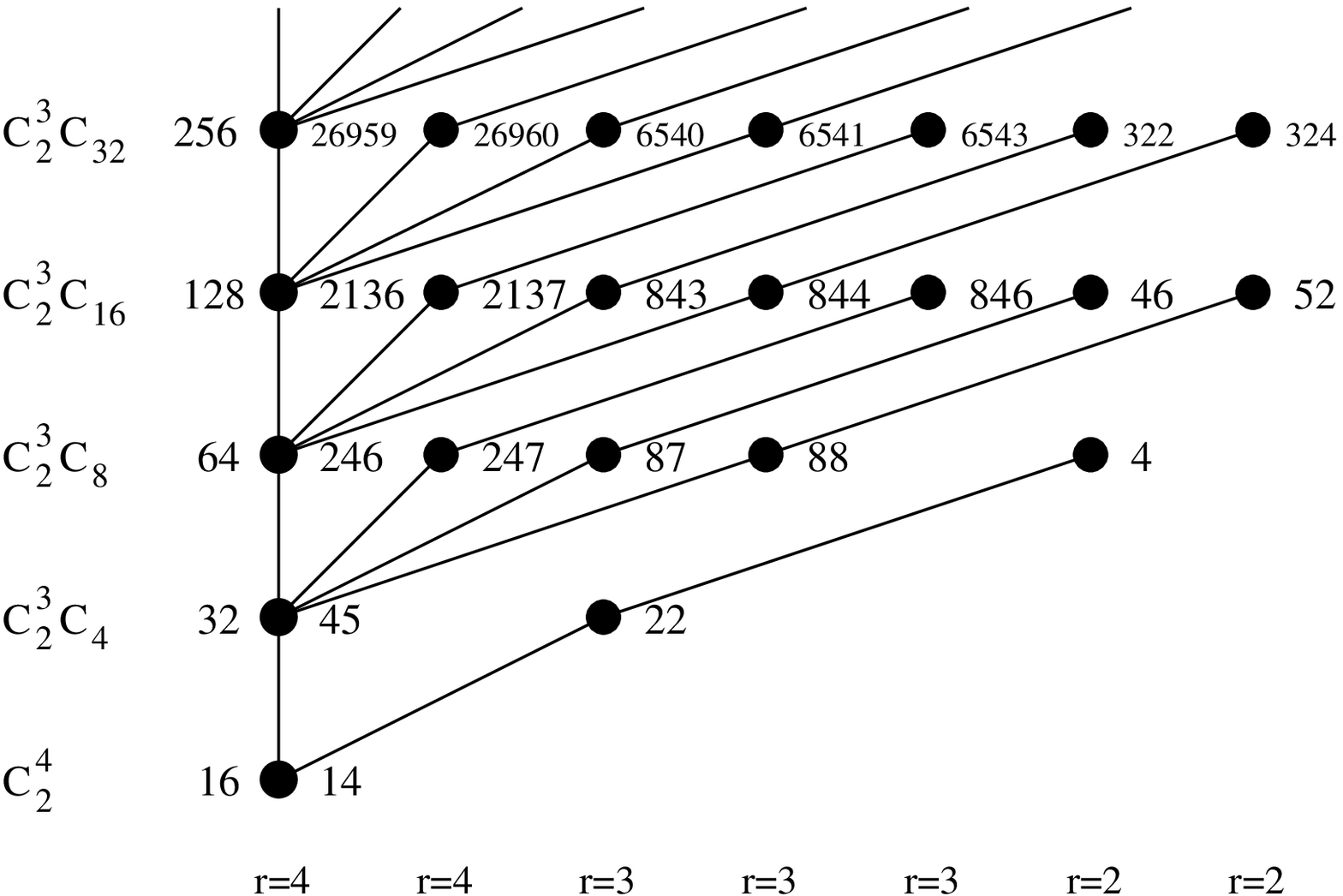} 
\hfill
\includegraphics[height=42mm]{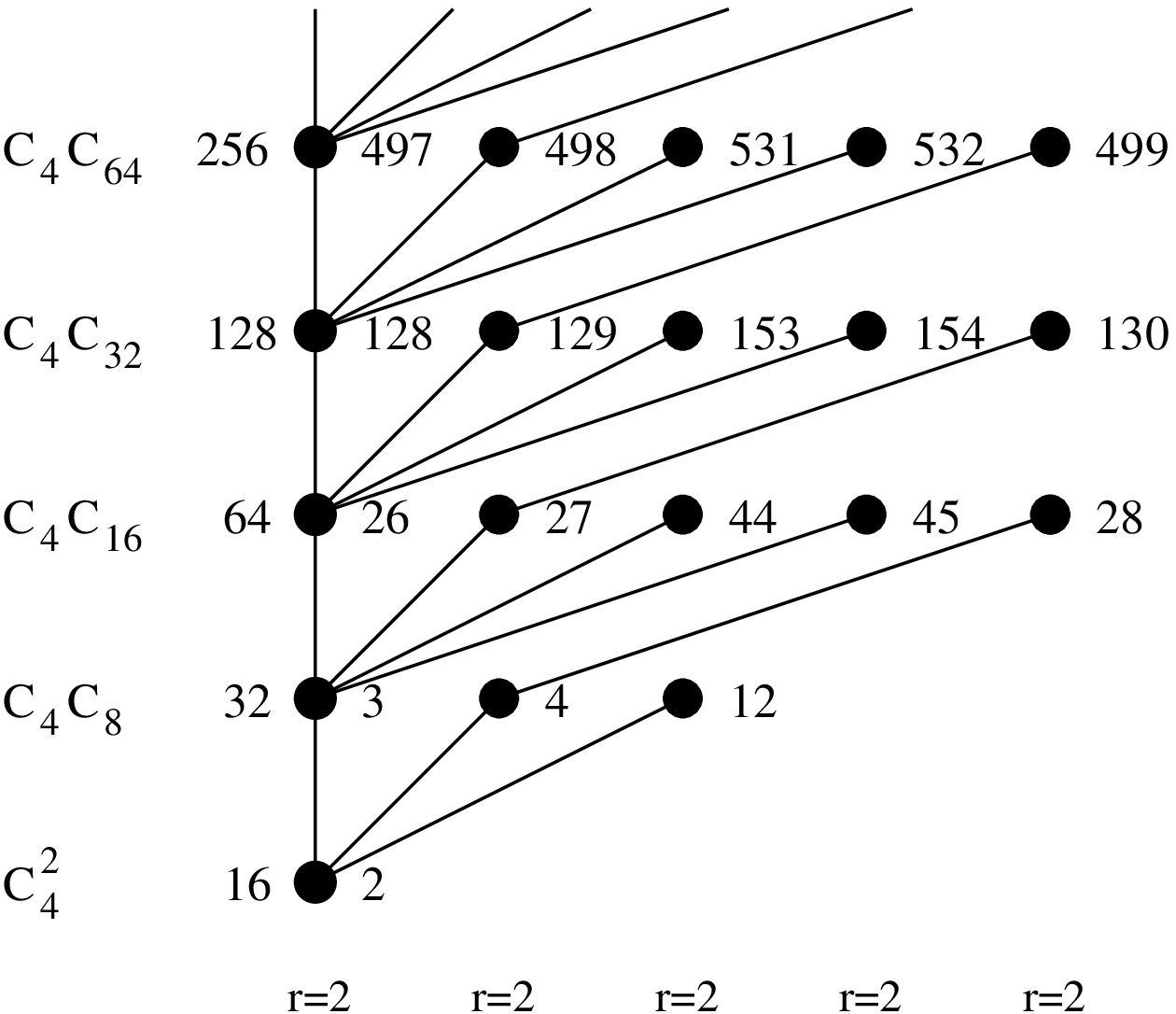} 

\abs\hrulefill
\end{table}

\AbsT{The trees rooted at non-abelian groups of order $16$.}
\label{2grpexmplcont}
\bfa
We have
$$ R:=\textsf{SmallGroups}(16,13)=(\lr{z}\tm\lr{y})\cn\lr{x} 
   \cong (C_4\tm C_2)\cn C_2 ,$$
where $Z(R)=\lr{z}$ and $R^{(1)}=\Phi(R)=\lr{z^2}$, such that $[y,x]=z^2$.
Hence we have $\exp(R)=\exp(Z(R))=4$ and $r(R)=3$, where 
$\cK(R)=\{w\in R;|w|\leq 2\}$ has precisely eight elements,
but is not a subgroup, thus $R$ is not of GK type.

\abs
Moreover, $\cZ(R)=\{z,z^3\}$ has precisely two elements, 
which are powers of each other, thus $\cT(R)$ has a unique stem.
Hence by \ref{stemthm} the unique trivial GK extension $G_l$ of $R$ 
of level $l\in\N$, thus having order $|G_l|=2^{l+4}$, is given as
$$ G_l=(\lr{g}\tm\lr{y})\cn\lr{x}\cong (C_{2^{l+2}}\tm C_2)\cn C_2,
\quad\text{where }Z(G_l)=\lr{g}\text{ and }g^{2^l}=z ;$$
moreover, since $\cZ(R)\cap\Phi(R)=\emp$ we from \ref{trivext}
get $r(G_l)=r(R)=3$.

\abs
Note that both the abelian group
$C_{2^{l+2}}\tm C_2\cong\lr{g}\tm\lr{y}\unlhd G_l$ and
the twisted dihedral group
$D_{2^{l+3}}^\ast\cong\lr{gy}\cn\lr{x}\unlhd G_l$
are maximal subgroups again of GK type, see
\ref{abelian} and \ref{cp2tree}, respectively,
even for $l=1$, while by construction we have
$\cK(G_l)=G_{l-1}=(\lr{g^2}\tm\lr{y})\cn\lr{x}\unlhd G_l$,
where we let $G_0:=R$.

\abs\bfb
We have 
$$ R':=\textsf{SmallGroups}(16,11)=(\lr{z'}\cn\lr{x'})\tm\lr{y'}
   \cong D_8\tm C_2 ,$$
where $|z'|=4$ and $|x'|=|y'|=2$, such that $[z',x']=(z')^2$.
Hence we have $Z(R')=\lr{(z')^2}\tm\lr{y'}$ and 
$(R')^{(1)}=\lr{(z')^2}$.
Moreover, we have
$$ R'':=\textsf{SmallGroups}(16,12)
       =\left(\lr{z''}\tm_{\lr{((z'')^2,(x'')^2)}}\lr{x''}\right)\tm\lr{y''} 
   \cong Q_8\tm C_2 ,$$
where $|z''|=|x''|=4$ and $|y''|=2$, such that $[z'',x'']=(z'')^2$.
Hence we have $Z(R'')=\lr{(z'')^2}\tm\lr{y''}$ and
$(R'')^{(1)}=\lr{(z'')^2}$.

\abs
By the comments on the groups $D_8$ and $Q_8$ in \ref{2grpexmpl}
we infer that both $R'$ and $R''$ are not of GK type either. 
Moreover, $\exp(R')=\exp(R'')=4$ and $\exp(Z(R'))=\exp(Z(R''))=2$ 
shows that the associated GK trees are finite. 
Indeed, it turns out that both trees have precisely two vertices.

\abs
We remark that, since the groups $D_8$ and $Q_8$ are isoclinic,
this also holds for the pair $R'$ and $R''$. We now show that
$R$ is isoclinic to $R'$ and $R''$ as well, 
where we only deal with the pair $R$ and $R'$, 
the argument for $R$ and $R''$ being analogous:
We have $R/Z(R)=\lr{\ov{y},\ov{x}}\cong C_2^2$,
where $\baraut\cn R\ra R/Z(R)$ denotes the natural epimorphism,
and $R'/Z(R')=\lr{\ov{z'},\ov{x'}}\cong C_2^2$,
where again $\baraut\cn R'\ra R'/Z(R')$ denotes the natural epimorphism.
Hence isoclinism is afforded by letting
$$ R^{(1)}\ra (R')^{(1)}\cn z^2\mt(z')^2 
\quad\text{and}\quad 
R/Z(R)\ra R'/Z(R')\cn\ov{y}\mt\ov{z'},\,\ov{x}\mt\ov{x'} .$$

\abs
Hence these examples show that isoclinic root groups might lead to
drastically different GK trees. Recall that by \ref{trivextdef}
all the groups lying on the stems of a fixed GK tree are mutually 
isoclinic, showing that going over to an isoclinic group does not 
necessarily preserve the property of being a root group either.

\begin{table}\caption{Trees rooted at non-abelian groups of order $16$.} 
\label{treetbl3}

\abs
\includegraphics[height=45mm]{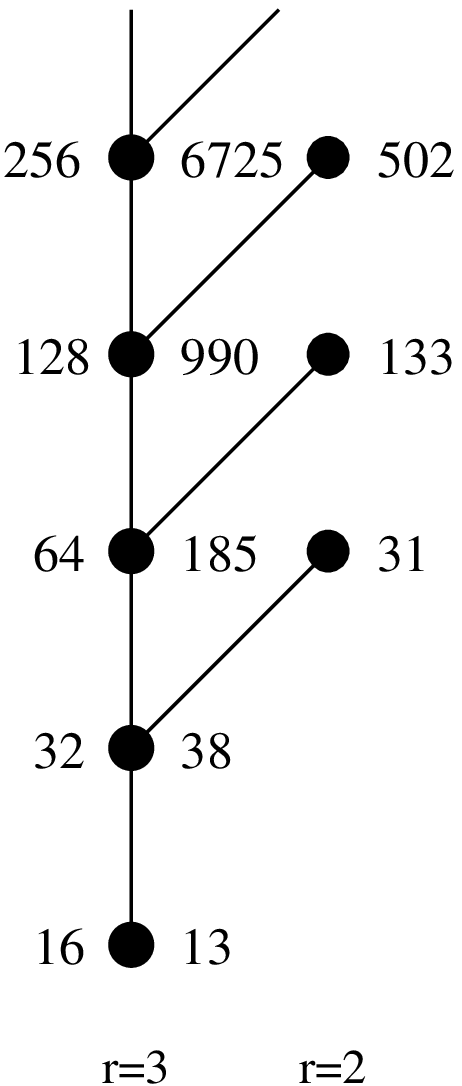} 
\hfill
\includegraphics[height=15mm]{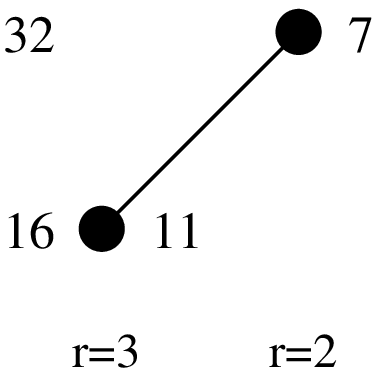} 
\hfill
\includegraphics[height=15mm]{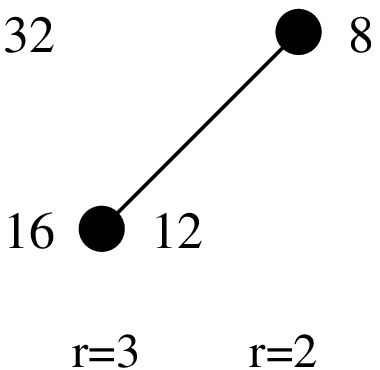} 
\hfill

\abs\hrulefill
\end{table}

\AbsT{Trees with multiple stems.}\label{stem2ex}
The smallest $2$-group being the root of
an infinite GK tree having more than one stem turns out to be 
$$ R:=\textsf{SmallGroups}(64,198)=\lr{w}\tm((\lr{z}\tm\lr{y})\cn\lr{x})
    \cong C_4\tm ((C_4\tm C_2)\cn C_2) ,$$
where the right hand direct factor is isomorphic to 
$\textsf{SmallGroups}(16,13)$, see \ref{2grpexmplcont}.

\abs
We point out that there is a similarity of this
example to the one given in \ref{stemexmpl} below,
which is obtained by replacing the direct factors $C_4$ and 
$(C_4\tm C_2)\cn C_2$ by $C_p$ and the extraspecial group 
$E_+(p^{2+1})$ of order $p^3$ and exponent $p$, respectively;
indeed this formal similarity is reflected in the structural analysis of 
both examples, but there are subtle differences.
We note that, just as $C_p\tm E_+(p^{2+1})$ can 
be generalised yielding the infinite series described in \ref{branchexmpl},
the present example $R$ is merely the first of a whole series,
but we will not delve into these constructions here, 
and just restrict ourselves to $R$:

\abs
By the comments in \ref{2grpexmplcont} on $\textsf{SmallGroups}(16,13)$,
which is not of GK type,
we have $\exp(R)=\exp(Z(R))=4$ and $r(R)=4$, where 
$Z(R)=\lr{w}\tm\lr{z}\cong C_4\tm C_4$ and 
$\Phi(R)=Z^2(R)=\lr{w^2}\tm\lr{z^2}$. Hence
we infer that $R$ is not of GK type either,
thus is the root of its infinite tree $\cT(R)$.

\abs
The set $\cZ(R)=Z(R)\smin Z^2(R)$ has $12$ elements. 
It turns out that $|\Aut(R)|=3072=2^{10}\cdot 3$, thus
$|\Out(R)|=768=2^8\cdot 3$, where $\Aut(R)$ acts on $\cZ(R)$ by 
a permutation group of order $32$, having two orbits:
$$ \begin{array}{ll}
\cO_1:=\{(w^{2i},z^j)\in\cZ(R);i\in\Z_2,j\in\Z_4^\ast\}, &
\text{with cardinality }4, \\
\cO_2:=\{(w^i,z^j)\in\cZ(R);i\in\Z_4^\ast,j\in\Z_4\}, &
\text{with cardinality }8. \rule{0em}{1.2em} \\
\end{array} $$
Moreover, both $\cO_1$ and $\cO_2$ are invariant under the
exponentiation action of $\Z_4^\ast$, as well as 
under the translation action of $Z^2(R)=\lr{w^2}\tm\lr{z^2}$,
hence coincide with the $A_l(R)$-orbits
on $\cZ(R)$, for all $l\geq 1$. 

\abs
Thus the tree $\cT(R)$ has two stems,
where branching of stems occurs at level $l=0$, that is
both stems just emanate from the root $R$.
Moreover, since both orbits we have
$\cO_1\cap\Phi(R)=\emp=\cO_2\cap\Phi(R)$,
by \ref{trivext} we have $r(G)=r(R)=4$ for all groups $G$
lying on a stem of $\cT(R)$.

\abs
The associated tree, as is found using \GAP{} and the
{\sf SmallGroups} database, is depicted in Table \ref{treetbl11};
here, we draw vertices as filled or open circles, depending on
whether the associated group has rank $4$ or $3$, respectively.

\begin{table}\caption{Tree rooted at $C_4\tm ((C_4\tm C_2)\cn C_2)$.} 
\label{treetbl11}

\rule{0mm}{179mm}
\hspace*{85mm}
\begin{rotate}{90}
\includegraphics[height=49mm]{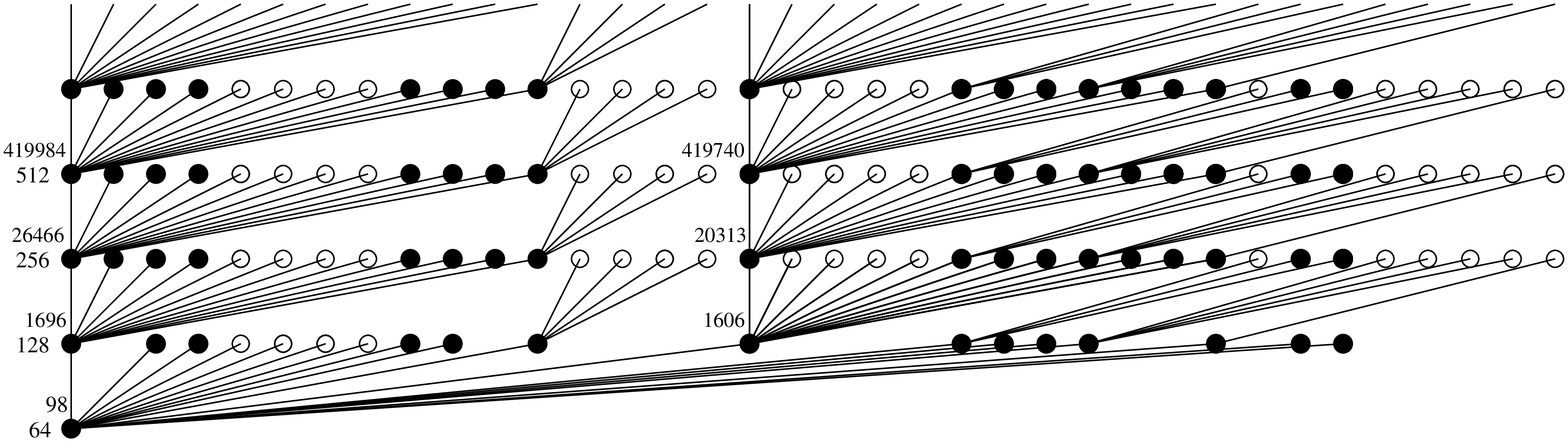} 
\end{rotate}

\abs\hrulefill
\end{table}

%% file: genex.tex
\section{Generic examples}\label{genex}

\abs
We now proceed towards generic GK trees, where the prime $p$ is 
treated as a parameter, and the root group is given as an abstract 
isomorphism type. Typically, the case $p=2$ needs special treatment
or has to be excluded.

\AbsT{The tree rooted at $C_p^2$.}\label{cp2tree}
We consider the elementary abelian group $C_p^2$, 
where $p$ is arbitrary, which is not of GK type.
Then, by \ref{abelianII} and \ref{trivextcor},
the unique stem of the tree $\cT(C_p^2)$ is
occupied on level $l\in\N$ by the abelian group 
$C_p\tm C_{p^{l+1}}$ and consists of the directed edges 
$C_p\tm C_{p^{l+1}}\ra C_p\tm C_{p^l}$. To determine the groups
not on lying on the stem, that is, the non-abelian groups in 
$\cT(C_p^2)$, we proceed as follows:

\abs
By \ref{immediate}, any group $G$ of GK type in $\cT(C_p^2)$ 
has cyclic deficiency $\dt(G)=\dt(C_p^2)=1$,
that is, if $G$ has level $l\in\N$, then it has order $|G|=p^{l+2}$ and
a cyclic maximal subgroup of order $\exp(G)=p^{l+1}$. The 
non-abelian $p$-groups of cyclic deficiency $1$ are well-known, 
see for example \cite[Thm.I.14.9]{hup}:

\abs\bfa
If $(p,l)\neq (2,1)$,
then let $G_l$ be the group of order $p^{l+2}$ given as
$$ G_l:=\lr{y}\cn\lr{g}\cong C_{p^{l+1}}\cn C_p,
\quad\text{where }y^g=y^{1+p^l} .$$
Hence $G_l$ has cyclic deficiency $1$. We show 
that $G_l$ is of GK type:

\abs
From $(y^p)^g=y^{p(1+p^l)}=y^p$ we deduce that
$H:=\lr{g,y^p}\cong C_p\tm C_{p^l}$ is a maximal subgroup
of $G_l$, having exponent $\exp(H)=p^l$. Moreover, for any 
$g^iy^j\in G_l$, where $i\in\Z_p$ and $j\in\Z_{p^{l+1}}$, we have 
$$ (g^iy^j)^p=g^{ip}y^{js}=y^{js},\quad\text{where }
   s:=\sum_{k=0}^{p-1}(1+p^l)^{ik} .$$
Thus for any prime $p$ and $l\geq 2$ we have 
$$ s\equiv\sum_{k=0}^{p-1}1=p\pmod{p^2} ,$$
while for $p$ odd and $l=1$ we have
$$ s\equiv\sum_{k=0}^{p-1}(1+ikp)=p+\binom{p}{2}\cdot ip
    \equiv p\pmod{p^2}. $$
Hence whenever $g^iy^j\in G_l\smin H$, that is
$j\in\Z_{p^{l+1}}^\ast$, we have 
$(g^iy^j)^p=y^{js}\in\lr{y^p}\smin\lr{y^{p^2}}$, implying that
$g^iy^j$ has order $p^{l+1}$. Thus we have $\cK(G_l)=H$.
This also shows that $G_l$ belongs to the GK tree containing
$H\cong C_p\tm C_{p^l}$, which is $\cT(C_p^2)$, and branches off
its stem at level $l-1$.

\abs
Now, if $p$ is odd, then 
the group $G_l$, where $l\geq 1$, is the only non-abelian group 
(up to isomorphism) of order $p^{l+2}$ of cyclic deficiency $1$,
implying that in this case all non-abelian groups of 
cyclic deficiency $1$ are of GK type. 

\abs
Moreover, this completes
the GK tree $\cT(C_p^2)$, see Table \ref{treegentbl1}: 
Next to its stem, it has bushes 
branching off at level $l-1$ for all $l\geq 1$,
consisting of a single non-stem vertex being connected to the stem
by the directed edge $C_{p^{l+1}}\cn C_p\ra C_p\tm C_{p^l}$,
in particular $\cT(C_p^2)$ is periodic from level $l=0$ on.
Finally, we point out that all the non-abelian
groups $G_l$, for $l\geq 1$, also have rank $r(G)=2=r(C_p\tm C_{p^l})$.
 
\abs\bfb
Hence it remains to consider the case $p=2$: 
Recall that for $l=1$ we only have the dihedral group $D_8$ and
the quaternion group $Q_8$, which have already been dealt with
in \ref{2grpexmpl}. Hence we may assume that $l\geq 2$. 
Then there are precisely four 
(isomorphism types of) non-abelian groups of order $2^{l+2}$ 
and cyclic deficiency $1$. They are given as
$$ G_{(x,a)}:=\lr{g,y\spmid g^2=x, y^{2^{l+1}}=1, y^g=y^a} ,$$
where $(x,a)\in G\tm\Z$ runs through the following cases:
$$ \begin{array}{|l|ll|}
\hline
(x,a) & G_{(x,a)} & \\
\hline
(1,-1) & D_{2^{l+2}}& \text{dihedral} \\
(1,-1+2^l) & SD_{2^{l+2}} & \text{semi-dihedral} \\
(1,1+2^l) & D_{2^{l+2}}^\ast & \text{twisted dihedral} \\
(y^{2^l},-1) & Q_{2^{l+2}}& \text{generalised quaternion} \\
\hline 
\end{array} $$

\abs
Note that the group $G_l$ above for the case $p=2$ and $l\geq 2$
yields the twisted dihedral group $D_{2^{l+2}}^\ast$,
which hence has already been shown to be of GK type,
being connected to the stem of $\cT(C_2^2)$ by the directed edge 
$D_{2^{l+2}}^\ast\ra C_2\tm C_{2^l}$.

\abs
We show that the remaining of the above groups are not of GK type:
Let $G_{(x,a)}$ be a dihedral, semi-dihedral or generalised quaternion 
group of order $2^{l+2}$, and hence of exponent $\exp(G_{(x,a)})=2^{l+1}$,
and assume that $G_{(x,a)}=\lr{g,y}$ is of GK type. Then we have $g^4=1$,
and from 
$$ (gy)^2=\left\{ \begin{array}{rl}
g^2=1, & \text{if }G_{(x,a)}\cong D_{2^{l+2}}, \\
y^{2^l}, & \text{if }G_{(x,a)}\cong SD_{2^{l+2}}, \\
g^2=y^{2^l}, & \text{if }G_{(x,a)}\cong Q_{2^{l+2}} \\
\end{array}\right. $$
we infer that $(gy)^4=1$ as well, hence since $l\geq 2$ we
get the contradiction
$$ G_{(x,a)}=\lr{g,gy}\leq\cK(G_{(x,a)})=\{w\in G_{(x,a)};|w|\leq 2^l\} .$$
Thus for $p=2$ the twisted dihedral group $D_{2^{l+2}}^\ast$
is the only non-abelian group of GK type (up to isomorphism) 
of order $2^{l+2}$ and cyclic deficiency $1$. 

\abs
Moreover, this completes
the GK tree $\cT(C_2^2)$, see Table \ref{treetbl1}: 
Next to its stem, it has bushes 
branching off at level $l-1$ for all $l\geq 2$,
consisting of a single non-stem vertex being connected to the stem
by the directed edge $D_{2^{l+2}}^\ast\ra C_2\tm C_{2^l}$,
while it has only one vertex in level $l=1$.
In particular, $\cT(C_2^2)$ is periodic from level $l=1$ on, 
and for all levels $l\geq 2$ coincides with the tree $\cT(C_p^2)$, 
for $p$ odd. 
Finally, we point out that all the non-abelian groups $D_{2^{l+2}}^\ast$,
for $l\geq 2$, have rank $r(G)=2=r(C_2^2)=r(C_2\tm C_{2^l})$.

\begin{table}\caption{Trees rooted at $C_p^2$ for $p$ odd.}
\label{treegentbl1}

\abs
\mbox{}\hfill
\includegraphics[height=60mm]{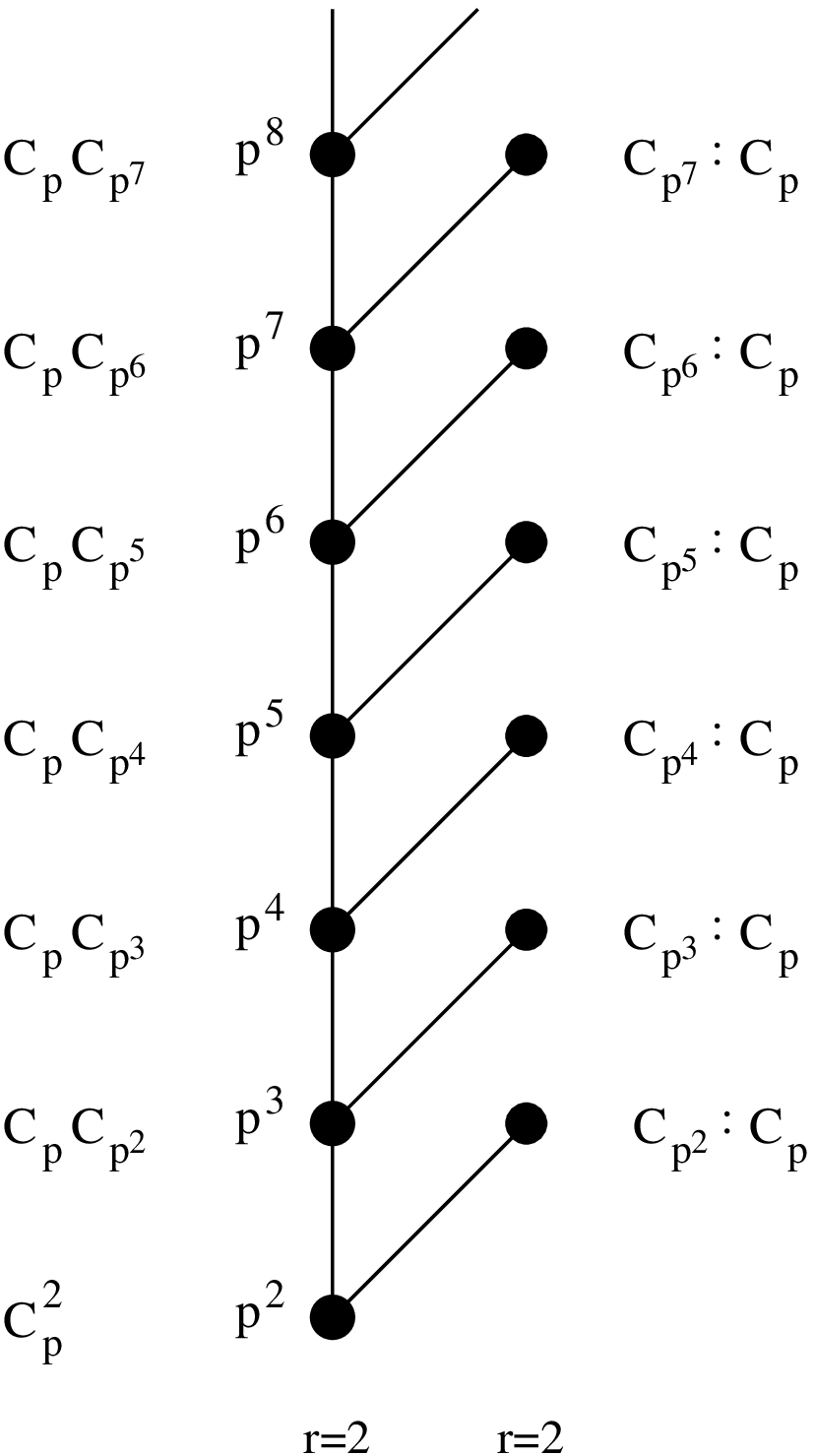} 
\hfill\mbox{}

\abs\hrulefill
\end{table}

\AbsT{The tree rooted at $C_p^3$.}\label{cp3tree}
Next we briefly consider the elementary abelian group 
$R:=\lr{x}\tm\lr{y}\tm\lr{z}\cong C_p^3$, where again $p$ is arbitrary,
which is not of GK type. Then, by \ref{abelianII} and \ref{trivextcor},
$\cT(R)$ has a unique stem, being occupied on level $l\in\N$ 
by the abelian group $C_p^2\tm C_{p^{l+1}}$, and consisting
of the directed edges $C_p^2\tm C_{p^{l+1}}\ra C_p^2\tm C_{p^l}$. 

\abs
At this stage, not having a theory describing the bushes in our hands,
we do not have to say much about the groups in $\cT(R)$
not on lying on the stem; recall that we have indicated
how $\cT(C_2^3)$ looks like in Table \ref{treetbl1}.
Here, we are content with exhibiting 
a non-abelian group of level $l=1$ in $\cT(R)$, in order to provide
generic examples of groups of GK type which are not powerful:

\abs
Let $G\cong R.C_p$ be the non-split non-trivial extension of degree $p$ 
of $R$ given as
$$ G:=\lr{x,y}\cn\lr{g}\cong C_p^2\cn C_{p^2},
\quad\text{where }[x,g]=y\text{ and }[y,g]=1, $$
and $R$ is naturally embedded into $G$ by letting $z:=g^p$.
Note that we hence have $r(G)=2$ and $Z(G)=\lr{y,z}\cong C_p^2$, 
in particular $\exp(Z(G))=p$.
Moreover, for all $i\in\Z_{p^2}$ and $j,k\in\Z_p$ we get
$$ (g^ix^jy^k)^p = z^i x^{jp} y^{kp+s}= z^i y^s \in R,
   \quad\text{where }s:=ij\cdot\frac{p(p-1)}{2} .$$
Hence from $\exp(R)=p$ we infer $\exp(G)=p^2$, and we have
$|g^ix^jy^k|=p^2$ if and only if $i\in\Z_{p^2}^\ast$,
or equivalently $g^ix^jy^k\in G\smin R$.
Hence $G$ is of GK type with kernel $\cK(G)=R$.

\abs
For $p=2$ we just recover $\textsf{SmallGroups}(16,3)$
in Table \ref{treetbl1}, where $G^4=\{1\}$ implies 
that $G$ is not powerful.
If $p$ is odd, then we have 
$G^p=\lr{z}\cong C_p$ and $G^{(1)}=\lr{y}\cong C_p$,
where $G^p\cap G^{(1)}=\{1\}$,
implying that $G$ is not powerful either.
\QED

\AbsT{Trees with stem branching.}\label{branchexmpl}
We now present examples of infinite GK
trees exhibiting interesting branching behaviour of their stems.
They are `doubly-generic' in the sense that we will
consider root groups of exponent $p^e$, where both the 
rational prime $p$ and $e\in\N$ are treated as parameters.

\abs
To begin with, let $p$ be odd, and let 
$$ E:=(\lr{z}\tm\lr{y})\cn\lr{x}\cong (C_{p^e}\tm C_{p^e})\cn C_{p^e},
\quad\text{where }[y,x]=z\text{ and }[z,x]=1 .$$
Hence we have $r(E)=2$ and $Z(E)=E^{(1)}=\lr{z}\cong C_{p^e}$. 
Moreover, from $y^x=yz$ we get $(y^j)^{(x^i)}=y^j z^{ij}$,
thus for all $i,j,k\in\Z_{p^2}$ and $a\in\N_0$ we have
$$ (x^i y^j z^k)^{p^a} = x^{ip^a} y^{jp^a} z^{kp^a+s(a)},
\quad\text{where }s(a):=ij\cdot\frac{p^a(p^a-1)}{2} .$$
Hence we get 
$E^{p^a}=(\lr{z^{p^a}}\tm\lr{y^{p^a}})\cn\lr{x^{p^a}}$;
note that $E^{p^a}$ is abelian if and only if $a\geq\cl{\frac{e}{2}}$.
Moreover, we have $\exp(E)=\exp(Z(E))=p^e$,
and $|x^i y^j z^k|<p^e$ if and only if $p\spmid\gcd(i,j,k)$, 
in other words if and only if $x^i y^j z^k\in E^p$; 
in particular $E$ is not of GK type.

\abs
Now let 
$$ R:=\lr{w}\tm((\lr{z}\tm\lr{y})\cn\lr{x})\cong C_{p^e}\tm E .$$
Hence we have $r(R)=3$ and $Z(R)=\lr{w}\tm\lr{z}\cong C_{p^e}^2$ 
and $R^{(1)}=E^{(1)}=\lr{z}\cong C_{p^e}$. Moreover, for $a\in\N_0$ we have
$R^{p^a}=\lr{w^{p^a}}\tm E^{p^a}$
and thus
$$ R^{p^a}R^{(1)}=\lr{w^{p^a}}\tm((\lr{z}\tm\lr{y^{p^a}})\cn\lr{x^{p^a}}) ;$$
in particular,
by Burnside's Basis Theorem \cite[Thm.III.3.15]{hup}, we have
$$ \Phi(R)=R^p R^{(1)}=
   \lr{w^p}\tm((\lr{z}\tm\lr{y^p})\cn\lr{x^p}) .$$
Since $E$ is not of GK type, neither is $R$,
and since we still have $\exp(R)=\exp(Z(R))=p^e$, 
we conclude that $R$ is the root of its infinite tree $\cT(R)$.

\abs
In order to describe the stems of $\cT(R)$,
we consider the action of $A_l(R)$ on $\cZ(R)$, for $l\in\{1\ld e\}$:
We have $|\cZ(R)|=p^{2e}-p^{2(e-1)}=p^{2e-2}(p^2-1)$, where
$$ \cZ(R)=Z(R)\smin Z^p(R)
=\{w^i z^j\in Z(R);i,j\in\Z_{p^e},\, p\spnmid\gcd(i,j)\} .$$

\abs
We first consider the exponentiation action of $\Z_{p^e}^\ast$:
A set of representatives of the $\Z_{p^e}^\ast$-orbits in $\cZ(R)$
is found by picking 
a generator of each of the cyclic subgroups of order $p^e$ in 
$Z(R)=\lr{w}\tm\lr{z}$. Hence considering their images with respect
to the projection map onto the left hand direct factor $\lr{w}$ yields
$$ \{w z^j\in\cZ(R);j\in\Z_{p^e}\}\dcup\coprod_{a\in\{1\ld e\}}
   \{w^{ip^a} z\in\cZ(R);i\in\Z_{p^{e-a}}^\ast\} .$$ 

\abs
Next, the group $R$, being is a trivial cyclic extension of $E$,
is a parent group in the sense of \ref{parents}. Hence we may consider 
the group $\Aut_0(R)\cong(\Aut(H)\tm\Z_{p^e}^\ast)\ltimes Z(H)$ 
of automorphism of $R$ leaving $H$ invariant, see \ref{automdirprod}:
Let $\psi\in\Aut_0(R)$ be the automorphism described by the triple
$(\id_H,1,z)\in(\Aut(H)\tm\Z_{p^e}^\ast)\ltimes Z(H)$,
and for all $k\in\Z_{p^e}^\ast$ let $\ph_k\in\Aut_0(R)$ 
be the automorphism described by the triple
$(\id_H,k,1)\in(\Aut(H)\tm\Z_{p^e}^\ast)\ltimes Z(H)$.

\abs
Then we have $(w z^j)^\psi=w z^{j+1}$, for all $j\in\Z_{p^e}$,
showing that the set $\{w z^j\in\cZ(R);j\in\Z_{p^e}\}$ is contained 
completely in a single $\Aut_0(R)$-orbit. Similarly, for $a\in\{1\ld e-1\}$
we have $(w^{ip^a} z)^{\ph_k}=w^{kip^a} z=w^{\ov{k}ip^a} z$, 
for all $i\in\Z_{p^{e-a}}^\ast$, where 
$\baraut\cn\Z_{p^e}\ra\Z_{p^{e-a}}$ denotes the natural epimorphism,
showing that the set $\{w^{ip^a} z\in\cZ(R);i\in\Z_{p^{e-a}}^\ast\}$
is contained completely in a single $\Aut_0(R)$-orbit as well.
Hence the set
$$ \{w\}\dcup\{w^{p^a} z;a\in\{1\ld e\}\} $$ 
contains a set of representatives of the 
$(\Aut_0(R)\tm\Z_{p^e}^\ast)$-orbits in $\cZ(R)$.

\abs
Now, for $a\in\N$ we have
$$ \cZ(R)\cap R^{p^a}R^{(1)}
  =\{w^{ip^a} z^j\in\cZ(R);i\in\Z_{p^{e-a}},j\in\Z_{p^e}^\ast\} .$$
Thus we conclude that the above elements belong to pairwise distinct 
$\Aut(R)$-orbits. Hence we have determined the 
$(\Aut(R)\tm\Z_{p^e}^\ast)$-orbits in $\cZ(R)$,
that is the $A_e(R)$-orbits. There are $e+1$ of them, 
given as follows, where $a\in\{1\ld e-1\}$:
$$ \begin{array}{lll}
\cO_0:=&\!\!\!\!\!\cZ(R)\smin R^pR^{(1)}
  &\!\!\!\!\!=\{w^i z^j\in\cZ(R);i\in\Z_{p^e}^\ast,j\in\Z_{p^e}\}, \\
\cO_a:=&\!\!\!\!\!\cZ(R)\cap\left(R^{p^a}R^{(1)}\smin R^{p^{a+1}}R^{(1)}\right) 
&\!\!\!\!\!=\{w^{ip^a} z^j\in\cZ(R);i\in\Z_{p^{e-a}}^\ast,j\in\Z_{p^e}^\ast\},
\rule{0em}{1.2em} \\
\cO_e:=&\!\!\!\!\!\cZ(R)\cap R^{p^e}R^{(1)}
  &\!\!\!\!\!=\{z^j\in\cZ(R);j\in\Z_{p^e}^\ast\}; \rule{0em}{1.2em} \\
\end{array} $$
their cardinalities are given as follows, where again $a\in\{1\ld e-1\}$:
$$ |\cO_0|=p^{2e-1}(p-1), \quad\quad
   |\cO_a|=p^{2e-2-a}(p-1),\quad\quad
   |\cO_e|=p^{e-1}(p-1)  .$$

\abs
Now we consider all the groups $A_l(R)$, where $l\in\{1\ld e\}$,
where we additionally have to take the translation action of 
$Z^{p^l}(R)=\lr{w^{p^l}}\tm\lr{z^{p^l}}$ on $\cZ(R)$ into account:
The orbit $\cO_0\sseq\cZ(R)$ is $Z^{p^l}(R)$-invariant, 
for all $l\in\{1\ld e\}$.
Moreover, for $i\in\{1\ld e\}$ the union
$$ \tcO_i:=\coprod_{a\in\{i\ld e\}}\cO_a\sseq\cZ(R) $$
is $Z^{p^l}(R)$-invariant, for all $l\in\{1\ld e\}$.
Hence, since $z\in\cO_e$ and $w^{p^a}z\in\cO_a$, for all $a\in\{1\ld e\}$,
we conclude that the 
$((\Aut(R)\tm\Z_{p^e}^\ast)\ltimes Z^{p^l}(R))$-orbits in $\cZ(R)$,
that is the $A_l(R)$-orbits, where $l\in\{1\ld e\}$, are given as
$$ \cO_0\dcup\cO_1\dcup\cdots\cO_{l-1}\dcup\tcO_l .$$

\abs
In conclusion, this shows that $\cT(R)$ has $e+1$ stems, being
parametrised by the $A_e(R)$-orbits $\cO_0\ld\cO_e\sseq\cZ(R)$.
More precisely, we have a two-fold branching of stems at level $l=0$,
described by the splitting $\cZ(R)=\cO_0\dcup\tcO_1$,
and we have further two-fold branching of stems at any level 
$l\in\{1\ld e-1\}$, described by the splitting 
$\tcO_l=\cO_l\dcup\tcO_{l+1}$.

\abs
Note that since
$\cO_0\cap\Phi(R)=\emp$ and $\cO_a\sseq\Phi(R)$, for all $a\in\{1\ld e\}$,
we infer from \ref{trivext} 
that all the GK groups lying on the stem belonging to $\cO_0$
have rank $r(R)=3$, while those 
lying on the stems belonging to $\cO_a$, where $a\in\{1\ld e\}$,
have rank $r(R)+1=4$.

\abs\abs
It seems to be worth-while to consider the smallest examples
of the above series more closely. In particular, they 
provide the first generic examples of non-abelian root groups,
and have cyclic deficiency $2$; recall that 
we have already covered the case of cyclic deficiency $1$ completely.

\AbsT{The case $e=1$.}\label{stemexmpl}
We keep the notation of \ref{branchexmpl}, in particular 
let still $p$ be odd, and let $e=1$.

\abs\bfa
Then we have $E=E_+(p^{2+1})$,
the extraspecial group of order $p^3$ and exponent $p$. 
As was already said, $E$ is not of GK type,
and thus is the root of its infinite tree $\cT(E)$.
Moreover, since $\Phi(E)=Z(E)=\lr{z}\cong C_p$ is cyclic of prime order, 
the set $\cZ(E)\sseq Z(E)$ consists 
of a single $A_1(E)$-orbit, hence $\cT(E)$ has a unique stem.
Since $\cZ(R)\sseq\Phi(E)$, by \ref{trivext} all GK groups $G$ lying on 
stem of $\cT(E)$ have rank $r(G)=r(R)+1=3$, and from \ref{immediate} 
we get $2\leq r(G)\leq 3$ for all groups $G$ in $\cT(R)$
which are not lying on the stem.

\abs
For example, for $p=3$ we have 
$E_+(3^{2+1})\cong\textsf{SmallGroups}(27,3)$,
and the associated tree, as is found using \GAP{} and the 
{\sf SmallGroups} database, is depicted in Table \ref{treetbl9};
note that the root has rank $r(E_+(3^{2+1}))=2$, which we
indicate by drawing the associated vertex as an open circle.

\abs
Actually, based on a few more experiments with \GAP{} and the
{\sf SmallGroups} database, we conjecture that the 
tree $\cT(E_+(p^{2+1}))$, for arbitrary odd $p$, 
coincides with $\cT(E_+(3^{2+1}))$,
at least for levels $l\geq 2$, and that all groups $G$ in
$\cT(E_+(p^{2+1}))$ not lying on the stem have rank $r(G)=2$.

\begin{table}\caption{Tree rooted at $E_+(3^{2+1})$.} 
\label{treetbl9}

\abs
\mbox{}\hfill
\includegraphics[height=45mm]{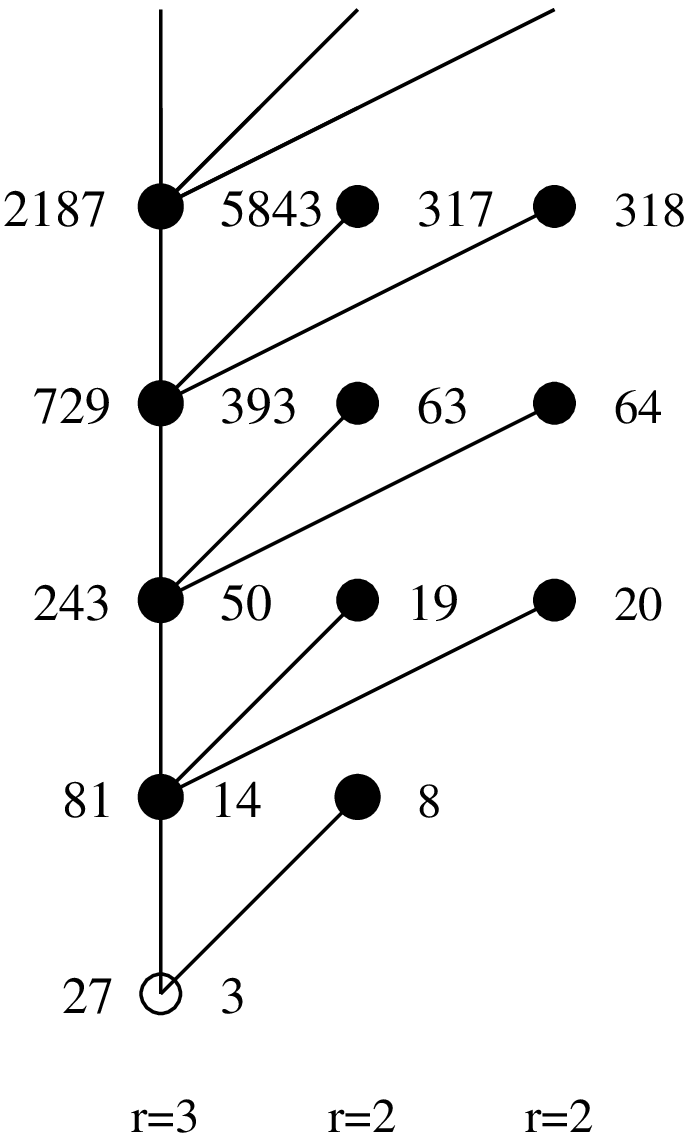} 
\hfill\mbox{}

\abs\hrulefill
\end{table}

\abs\bfb
We finally briefly consider $R=C_p\tm E$.
For example, for $p=3$ we have
$$ C_3\tm E_+(3^{2+1})\cong\textsf{SmallGroups}(81,12) ,$$
and the associated tree, as is found using \GAP{} and the 
{\sf SmallGroups} database, is depicted in Table \ref{treetbl10}.
Based on few more experiments with \GAP{} and the
{\sf SmallGroups} database, we conjecture that 
the trees $\cT(C_p\tm E_+(p^{2+1}))$, for arbitrary odd $p$,
are very similar to $\cT(C_3\tm E_+(3^{2+1}))$, but might differ
in details.

\begin{table}\caption{Tree rooted at $C_3\tm E_+(3^{2+1})$.} 
\label{treetbl10}

\abs
\mbox{}\hfill
\includegraphics[height=42mm]{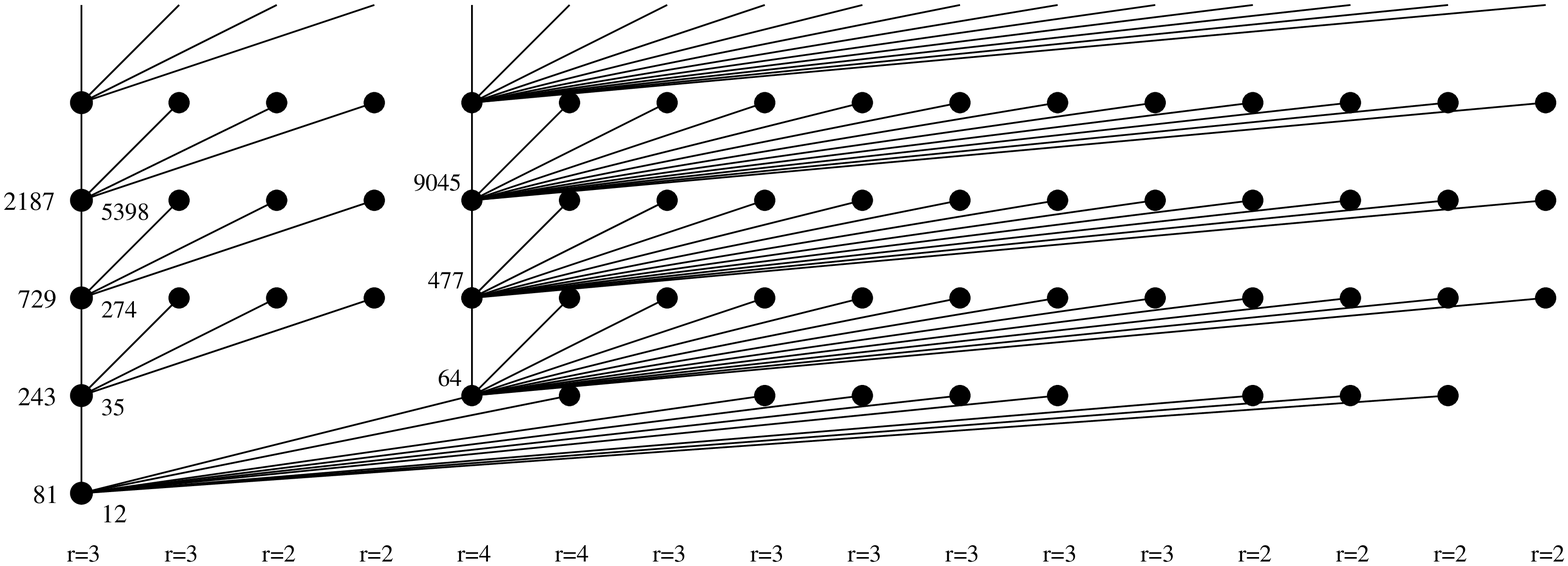} 
\hfill\mbox{}

\abs\hrulefill
\end{table}